\documentclass{scrartcl}
\usepackage{enumitem}
\usepackage{mathtools}

\usepackage{tikz}
\usetikzlibrary{decorations.pathreplacing}

\usepackage{amsmath}
\usepackage{amsfonts}
\usepackage{amssymb}
\usepackage{amsthm}
\usepackage{dsfont}
\usepackage[T1]{fontenc}
\usepackage[latin1]{inputenc}
\usepackage{epsfig}
\usepackage{graphicx}
\usepackage{bm}
\usepackage[english]{babel}

\usepackage[boxed, lined]{algorithm2e}
\usepackage{float}
\usepackage{comment}

\newcommand{\E}{\mathbf{E}}

\newcommand{\bk}{\mathbf{k}}
\newcommand{\bj}{\mathbf{j}}
\newcommand{\bi}{\mathbf{i}}

\newcommand{\bl}{\mathbf{l}}
\newcommand{\bx}{\mathbf{x}}

\newcommand{\bz}{\mathbf{z}}

\newcommand{\bv}{\mathbf{v}}

\newcommand{\bX}{\mathbf{X}}

\newcommand{\D}{{\mathcal{D}}}
\newcommand{\A}{{\mathcal{A}}}

\newcommand{\F}{{\mathcal{F}}}

\newcommand{\Nu}{{\mathcal{N}}}
\newcommand{\M}{{\mathcal{M}}}

\newcommand{\N}{\mathbb{N}}

\newcommand{\R}{\mathbb{R}}
\newcommand{\Z}{\mathbb{Z}}
\newcommand{\Rd}{\mathbb{R}^d}

\newcommand{\beq}{\begin{eqnarray*}}

\newcommand{\eeq}{\end{eqnarray*}}

\newcommand{\beqm}{\begin{eqnarray}}

\newcommand{\eeqm}{\end{eqnarray}}

\newtheorem{theorem}{Theorem}

\newtheorem{lemma}{Lemma}
\newtheorem{definition}{Definition}

\DeclareMathOperator*{\argmin}{arg\,min}

\newcommand{\EXP}{{\mathbf E}}
\newcommand{\PROB}{{\mathbf P}}

\renewcommand{\P}{{\cal P}}
\renewcommand{\bf}{\normalfont \bfseries}
\renewcommand{\it}{\normalfont \itshape}

\allowdisplaybreaks
\begin{document}
\begin{center}

  {\LARGE \bf
    Estimation of a regression function on a manifold by
    fully connected deep neural networks}\footnote{
Running title: {\it Estimation of a regression function on a manifold }}
\vspace{0.5cm}

Michael Kohler, Sophie Langer$\footnote{Corresponding author. Tel:
  +49-6151-16-23371}$ and Ulrich Reif\\

{\it 
Fachbereich Mathematik, Technische Universit\"at Darmstadt,
Schlossgartenstr. 7, 64289 Darmstadt, Germany,
email: kohler@mathematik.tu-darmstadt.de,
langer@mathematik.tu-darmstadt.de, reif@mathematik.tu-darmstadt.de }

\end{center}
\vspace{0.5cm}

\begin{center}
	April 11, 2021
\end{center}
\vspace{0.5cm}

\noindent
    {\bf Abstract}\\
    Estimation of a regression function from  independent and
    identically distributed data
    is considered. The $L_2$ error with integration with respect to
    the distribution of the predictor variable is used as the error
    criterion. The rate of convergence of least squares estimates
    based on fully connected spaces of deep neural networks
    with ReLU activation function is analyzed for smooth regression functions.
    It is shown that in case that the distribution of the
    predictor variable is concentrated on a manifold,
    these estimates achieve a rate of convergence which depends
    on the dimension of the manifold and not on the number
    of components of the predictor variable.
    \vspace*{0.2cm}


\vspace*{0.2cm}


\section{Introduction}
\label{se1}
Deep neural networks (DNNs) are built of multiple layers
and learn sequentially multiple levels of representation
and abstraction by performing a nonlinear transformation on the data. 
The approach has 
proven itself to work incredibly well in practice, like for speech 
(Graves et al.\ (2013)) and image recognition (Krizhevsky et al.\
(2017)), or game intelligence (Silver et al.\ (2016)).  
But, unfortunately, the procedure is not well understood.  Recently, 
several researchers tried to explain the performance of DNNs from a 
theoretical point of view.  Results concerning the approximation power 
of DNNs were shown in Montufar (2014), Eldan and Shamir (2016), 
Yarotsky (2017),  Yarotsky and Zhevnerchuck (2020),  Langer 
(2021b) and Lu et al.  (2020).  Beside this,
quite a few articles try to answer the question about why neural 
networks perform well on unknown new data sets
(cf., e.g., Bauer and Kohler (2019), Schmidt-Hieber (2020), 
Kohler and Langer (2020), Kohler, Krzy\.zak, and Langer (2019), 
Langer (2021a), Imaizumi and Fukumizu (2019), Suzuki (2018), Suzuki 
and Nitanda (2019), Schmidt-Hieber (2019) and the literature cited 
therein).
The standard framework to do this is
to consider DNNs in the context of nonparametric regression. Here,
$(\bX,Y)$ is an $\Rd \times \R$--valued random vector
satisfying 
$\EXP \{Y^2\}<\infty$, and given a sample
of size $n$ of $(\bX,Y)$, i.e., given a data set
\begin{equation*}
\D_n = \left\{
(\bX_1,Y_1), \ldots, (\bX_n,Y_n) 
\right\},
\end{equation*}
where
$(\bX,Y)$, $(\bX_1,Y_1)$, \ldots, $(\bX_n,Y_n)$ are
independent and
    identically distributed
(i.i.d.),
the aim is to construct an estimate
\[
m_n(\cdot)=m_n(\cdot, \D_n):\Rd \rightarrow \R
\]
of the so--called regression function $m:\Rd \rightarrow \R$,
$m(\bx)=\EXP\{Y|\bX=\bx\}$ such that the so--called $L_2$ error
\[
\int |m_n(\bx)-m(\bx)|^2 {\PROB}_{\bX} (d\bx)
\]
is ``small'' (cf., e.g., Gy\"orfi et al.\ (2002)
for a systematic introduction to nonparametric regression and
a motivation for the $L_2$ error).

It is well-known that without smoothness assumptions on the
regression function it is not possible to derive nontrivial results
on the rate of convergence of nonparametric regression
estimates (cf.\ Cover (1968) and Section~3 in Devroye and Wagner 
(1980)).
In the sequel we assume that the regression function is
$(p,C)$-smooth according to the following definition.
\begin{definition}
\label{intde1} 
  Let $p=q+s$ for some $q \in \N_0$ and $0< s \leq 1$.
A function $m:\R^d \rightarrow \R$ is called
{\em $(p,C)$-smooth}, if for every $\bm{\alpha}=(\alpha_1, \dots, 
\alpha_d) \in \N_0^d$
with $\sum_{j=1}^d \alpha_j = q$ the partial derivative
$\partial^q m/(\partial x_1^{\alpha_1}
\dots
\partial x_d^{\alpha_d}
)$
exists and satisfies
\[
\left|
\frac{
\partial^q m
}{
\partial x_1^{\alpha_1}
\cdots
\partial x_d^{\alpha_d}
}
(\bx)
-
\frac{
\partial^q m
}{
\partial x_1^{\alpha_1}
\cdots
\partial x_d^{\alpha_d}
}
(\bz)
\right|
\leq
C
\| \bx-\bz \|^s
\]
for all $\bx,\bz \in \R^d$, where $\Vert\cdot\Vert$ denotes the Euclidean norm.
\end{definition}
As it was shown by Stone (1982),
the optimal Minimax rate of convergence for estimation of a
$(p,C)$--smooth regression function is $n^{-2p/(2p+d)}$.  
This rate underlies one big problem, namely the so-called 
\textit{curse of dimensionality}: For fixed $p$ and increasing $d$ 
this rate gets rather slow.  Since $d$ tends to be very large in many 
machine learning applications, 
to show Stone's Minimax rate for DNNs is not really the answer for 
the empirical good performance of DNNs.
This is why many results
are restricted to further assumptions.  
Bauer and Kohler (2019), Schmidt-Hieber (2020), and Langer (2021a) 
considered regression functions with some kind of compositional 
structure and showed that DNNs achieve a dimensionality reduction
in this setting.
Kohler, Krzy\.zak, and Langer (2019),
and Eckle and Schmidt-Hieber (2019) could show that
DNNs can mimic the form of multivariate adaptive regression splines (MARS).
Kohler, Krzy\.zak, and Langer (2019) further showed that in case of 
regression functions with low local dimensionality DNNs are able to 
achieve dimensionality reduction.  Further approaches like Barron 
(1993,1994), Suzuki (2018), and Suzuki and Nitanda (2019) consider 
various types of smoothness or spectral distibutions.  

\subsection{Intrinsic Dimensionality}
All the above mentioned results mainly focus on the structure of the 
underlying regression function. Less results explore the geometric 
properties of the data.  But it is reasonable to also focus on the 
structure of the input variable $\bX$. Firstly, there exist several 
results where high-dimensional problems can be treated in much lower 
dimension
(cf., e.g., Tenenbaum et al.\ (2000)). For instance, if we consider 
the pixels of potraits of persons, the input dimensionality may be 
quite high, but the meaningful structure of these images and therefore 
the intrinsic dimensionality can lie in a much smaller space. 
Secondly, as already mentioned in Imaizumi and Nakada (2020), many 
estimators like kernel methods or Gaussian process regression show 
good rate of convergence results depending only on the intrinsic 
dimensionality of the input data (cf., e.g., Bickel and Li (2007) and 
Kpotufe (2011)).
It is therefore interesting to investigate whether estimators based on DNNs
are also able to exploit the structure of the input data.  
\\
\\
In the sequel we
do this by 
considering the special case that $\bX$ is 
concentrated on some $d^*$--dimensional Lipschitz-manifold. To 
describe this formally, we use the following definition.
\begin{definition}
  \label{intde2} 
  Let $\M \subseteq \Rd$ be compact and let $d^* \in \{1, \dots, d\}$.

  \noindent
      {\bf a)} We say that $U_1,\dots,U_r$ is
      an {\em open covering  of $\M$}, if $U_1,\dots,U_r \subset \Rd$
      are open (with respect to the Euclidean topology on $\Rd$)
      and satisfy
      \[
\M \subseteq \bigcup_{l=1}^r U_l.
\]

  \noindent
      {\bf b)} We say that
\[
\psi_1, \dots, \psi_r: [0,1]^{d^*}\rightarrow \Rd
\]
are {\em bi-Lipschitz functions}, if there exists $0 < C_{\psi,1} 
\leq C_{\psi,2}
< \infty$ such that
\begin{equation}
\label{de2eq1}
C_{\psi,1}  \cdot \|\bx_1-\bx_2\|
\leq
\| \psi_l(\bx_1)-\psi_l(\bx_2) \|
\leq
C_{\psi,2}  \cdot \|\bx_1-\bx_2\|
\end{equation}
holds for any $\bx_1,\bx_2 \in  [0,1]^{d^*}$ and any
$l \in \{1, \dots, r\}$.

        \noindent
            {\bf c)}
            We 
say that $\M$ is a {\em $d^*$-dimensional Lipschitz-manifold} if 
there exist 
bi-Lipschitz functions $\psi_i : [0,1]^{d^*} \to \R^d$ $(i \in \{
  1,\dots,r \})$, 
and an open covering $U_1,\dots,U_r$ of $\M$ such that 
\[
\psi_l( (0,1)^{d^*} ) = \M \cap U_l
\]
holds for all $i \in \{1, \dots, r\}$.
Here we call $\psi_1, \dots, \psi_r$ the {\em parametrizations}
of the manifold.
\end{definition}

If $\M$ is a $d^*$-dimensional Lipschitz-manifold, then
\begin{equation}
  \label{de2eq2}
\M = \bigcup_{l=1}^r \M \cap U_l
=
\bigcup_{l=1}^r \psi_l\bigl( (0,1)^{d^*} \bigr)
\end{equation}
and (\ref{de2eq1}) hold. We will see in the proof of our main result
that it suffices to assume that $\M$ has these two properties.

\subsection{Neural Networks}
In our analysis we consider DNNs with rectifier linear unit 
(ReLU)  $\sigma(x) = \max\{x,0\}$ as activation function.  A neural
networks can be described by its {\em number of hidden layers} 
$L \in \N$ and its {\em number of neurons per layer} $\textbf{k} = 
(k_1, \ldots, k_{L}) \in \mathbb{N}^{L}$, where $k_i$ describes the 
number of neurons in the $i$-th layer.  This leads to a network 
architecture $(L, \bk)$ and the corresponding neural network can be 
defined as follows:
\begin{definition}
A {\em multilayer feedforward neural network with network 
architecture $(L, \textbf{k})$ and ReLU activation function $\sigma$}
is a real-valued function defined on $\mathbb{R}^d$ of the form
\begin{equation}\label{inteq1}
f(\bx) = \sum_{i=1}^{k_L} c_{1,i}^{(L)}f_i^{(L)}(\bx) + c_{1,0}^{(L)}
\end{equation}
for some $c_{1,0}^{(L)}, \ldots, c_{1,k_L}^{(L)} \in \mathbb{R}$ and 
for functions $f_i^{(L)}$ recursively defined by
\begin{equation}
  \label{inteq2}
f_i^{(s)}(\bx) = \sigma\left(\sum_{j=1}^{k_{s-1}} c_{i,j}^{(s-1)} f_j^{(s-1)}(\bx) + c_{i,0}^{(s-1)} \right)
\end{equation}
for some $c_{i,0}^{(s-1)}, \dots, c_{i, k_{s-1}}^{(s-1)} \in \mathbb{R}$,
$s \in \{2, \dots, L\}$,
and
\begin{equation}
  \label{inteq3}
f_i^{(1)}(\bx) = \sigma \left(\sum_{j=1}^d c_{i,j}^{(0)} x^{(j)} + c_{i,0}^{(0)} \right)
\end{equation}
for some $c_{i,0}^{(0)}, \dots, c_{i,d}^{(0)} \in \mathbb{R}$. 
\end{definition}

%
%
%
%
%
The space of multilayer neural networks with 
$L$ hidden layers and $r$ neurons per layer 
is defined by
\begin{align}\label{F}
  \mathcal{F}(L, r) = \{ &f \, : \,  \text{$f$ is of the form } \eqref{inteq1}
  \text{ with }
k_1=k_2=\ldots=k_L=r 
\}.
\end{align}
As in Kohler and Langer (2020), we denote this network class as 
\textit{fully connected neural networks}.  In contrast, network 
classes
with a further restriction on the total number of nonzero weights 
in the network are called {\em sparsely connected neural networks}. 
\\

A corresponding least squares estimator can  be defined by 
\begin{align*}
\tilde{m}_n(\cdot) = \argmin_{f \in \mathcal{F}(L_n,r_n)} \frac{1}{n} 
\sum_{i=1}^n \left|f(\bX_i) - Y_i\right|^2.
\end{align*}
Here the number of hidden layers $L_n$ and the number of neurons $r_n$
is chosen in dependence of the sample size.  For simplicity we assume
here and in the sequel that the minimum above indeed exists.  When
this is not the case our theoretical results also hold for any
estimator which minimizes the above empirical $L_2$ risk up to a small additional term.
\\
Due to the fact that we do not impose any sparsity constraints (i.e., 
an additional bound on the number of nonzero weights in the functions
in our neural network space),  the least squares estimate above
can be (approximately) implemented in a very simple way with standard 
software
packages, e.g., in the Deep Learning framework of {\em tensorflow} 
and
{\em keras}. Here approximating the above least squares estimate
from data $x_{\mathrm learn}$ and $y_{\mathrm learn}$ can be done 
with only a few 
lines of code as follows:
\begin{quote}
 model = Sequential() \\
 model.add(Dense(d, activation="relu", input\_shape=(d,))) \\             
 for i in np.arange(L):\\
 \hspace*{3cm}     model.add(Dense(K, activation="relu")) \\
 model.add(Dense(1)) \\
 model.compile(optimizer="adam",
                  loss="mean\_squared\_error") \\
 model.fit(x=x\_learn,y=y\_learn)
  \end{quote}
For the implementation of sparsely connected neural networks,  one may 
use so-called {\em pruning methods},  
which start with large strongly connected neural networks and delete 
redundant parameters during the training process. Although the 
procedure is also simple, the computational costs of these methods are 
high, due to the large initial size of the networks. Therefore, the 
implementation of sparsely connected networks is critically questioned 
in the literature (see, e.g., Evci et al.\ (2019) or Liu et al.\
(2019)).

%
%

\subsection{Main results}
In our study,  we analyze the performance of DNNs when the input 
values are concentrated on a $d^*$-dimensional Lipschitz-manifold.
Here we develop a proof technique
which enables us to show a rate of convergence which depends
 only on the smoothness $p$ and dimension $d^*$ of our manifold, but not on $d$.
In particular, we show
that the convergence rate of the above least squares
neural network
estimator is of the order $n^{-2p/(2p+d^*)}$
(up to some logarithmic factor)
and therefore able to circumvent the curse of dimensionality in case 
that $d^*$ is rather small.
In contrast to earlier results (see, e.g., Schmidt-Hieber (2019) and 
Imaizumi and Nakada (2020)) 
we consider fully connected DNNs in our analysis.  In these networks  
the number of hidden layers is bounded by some logarithmic factor in 
the sample size and the number of neurons per layer tends to infinity 
suitably fast for sample size tending to infinity.
To show the above rate of convergence we derive a
new approximation result concerning fully connected DNNs 
on $d^*$ -dimensional Lipschitz-manifolds.
As in Kohler and Langer (2020) we use a two-step approximation, where 
we partition our space in a finite set of coarse and fine hypercubes 
and approximate the Taylor polynomial on each of the cubes by a DNNs.  
  As already mentioned
above the analysis of fully connected DNNs has the main advantage that 
an approximate implementation of a corresponding least squares 
estimate 
is much easier.
\subsection{Related results}
First results concerning neural networks date back to the 1980's.  
Here networks with only one hidden layer, so-called shallow neural 
networks, were analyzed and it was shown that they approximate any continuous function arbitrarily well provided the 
number of neurons is large enough (see, e.g., Cybenko (1989)).  Limits 
of this network architecture were analyzed e.g.\ in Mhaskar and 
Poggio (2016), who showed that specific functions cannot be 
approximated by shallow neural networks but by deep ones.
Concerning results on the approximation power of
multilayer neural networks we refer to
Montufar (2014), Eldan and Shamir (2016),
Yarotsky (2017),
Yarotsky and Zhevnerchuck (2019), Langer (2021b) and the literature cited therein.
The generalization error of least squares estimates based on sparsely connected
multilayer neural networks has been investigated in
Bauer and Kohler (2019),
Kohler and Krzy\.zak (2017),
Schmidt-Hieber (2020),  Imaizumi and Fukumizu (2020), and Suzuki 
(2018). In particular, it was shown that DNNs can achieve a dimension 
reduction in case that the
regression function is a composition of (sums of) functions,  where
each of the function depends only on a few variables.
As was shown in Kohler and Langer (2020) and Langer (2021a)
similar results can also be achieved for least squares estimates based
on fully-connected multilayer neural networks, which are easier to implement.
\\
\\
Function approximation and estimation on manifolds has been studied 
e.g.\ in Belkin and Niyogi (2008),  
Singer (2006),  Davydov and Schumaker (2007), Hangelbroek, Narcowich 
and Ward (2012), and Lehmann et al.\ (2019).  An analysis in 
connection to DNNs is given in Mhaskar (2010), who showed an 
approximation rate using so-called Eignets.  An overview of related 
results can be found in Chui and Mhaskar (2018).  In Schmidt-Hieber 
(2019) approximation rates and statistical risk bounds for functions 
defined on a manifold were derived.  His result is restricted to 
sparsely connected neural networks, where only a bounded number of 
parameters in the network is non-zero.    In comparison to our result 
his proof strategy is more complex and requires a stronger smoothness 
assumption on the
charts of the manifold. In Imaizumi and Nakada (2020) approximation 
rates and statistical risk bounds depending on the Minkowski dimension 
of the domain were shown, which is a more general framework than in 
our paper
since the Minkowski dimension of the $d^*$-dimensional Lipschitz-manifolds
considered in our paper is bounded from above by $d^*$.
In case that $\operatorname{supp}(\bX)$ is only of dimension $d^* < 
d$ in the Minkowski sense,
Imaizumi and Nakada (2020) derived similar to us a $d^*$-dimensional 
approximation rate and a convergence rate of a corresponding
least squares estimator of $n^{-2p/(2p+d^*)}$.
But in contrast to our results, 
the DNNs in Imaizumi and Nakada (2020) were restricted by a further 
sparsity constraint, such that an implementation of a corresponding 
estimator is more difficult.

In the computer science literature the problem
considered in this paper is called {\em manifold learning}, 
see Subsection~5.11.3 in
Goodfellow, Bengio, and Courville (2016) and the literature cited 
therein.
%
%
%
%

\subsection{Notation}
\label{se1sub7}
Throughout the paper, the following notation is used:
The sets of natural numbers, natural numbers including $0$,
integers, and real numbers
are denoted by $\N$, $\N_0$, $\Z$, and $\R$, respectively.
For $z \in \R$, we denote
the smallest integer greater than or equal to $z$ by
$\lceil z \rceil$, and $\lfloor z \rfloor$ denotes
the largest integer less than or equal to $z$.
Vectors are denoted by bold letters, e.g., $\bx = (x^{(1)}, \dots, 
x^{(d)})^T$. 
We define $\bold{1}=(1, \dots, 1)^T$ and $\bold{0} = (0, \dots, 0)^T$.
A $d$-dimensional multi-index is a $d$-dimensional vector $\bold{j} =
(j^{(1)}, \dots, j^{(d)})^T \in \N_0^d$. As usual, we define
$\|\bold{j}\|_1 = j^{(1)}+\dots+j^{(d)}$, $\bold{j}! = j^{(1)}! \cdots
j^{(d)}!$, and
\[
\partial^{\bold{j}} = \frac{\partial^{j^{(1)}}}{\partial (x^{(1)})^{j^{(1)}}} \cdots \frac{\partial^{j^{(d)}}}{\partial (x^{(d)})^{j^{(d)}}}.
\]
Let $D \subseteq \R^d$ and let $f:\R^d \rightarrow \R$ be a real-valued
function defined on $\R^d$.
We write $\bx = \argmin_{\bz \in D} f(\bz)$ if
$\min_{\bz \in \D} f(\bz)$ exists and if
$\bx$ satisfies
$\bx \in D$ and $f(\bx) = \min_{\bz \in \D} f(\bz)$.
For $f:\R^d \rightarrow \R$,
\[
\|f\|_\infty = \sup_{\bx \in \R^d} |f(\bx)|
\]
is its supremum norm, and the supremum norm of $f$
on a set $A \subseteq \R^d$ is denoted by
\[
\|f\|_{\infty, A} = \sup_{\bx \in A} |f(\bx)|.
\]
Furthermore, we define the norm $\| \cdot \|_{C^q(A)}$ of the smooth 
function space $C^q(A)$ by
\begin{align*}
\|f\|_{C^q(A)} =\max\left\{\|\partial^{\bj}f\|_{\infty, A}: \|\bj\|_1 
\leq q, \bj \in \N^d\right\}
\end{align*}
for any $f \in C^q(A)$.  Let $\F$ be a set of functions $f:\Rd \rightarrow \R$,
let $\bx_1, \dots, \bx_n \in \Rd$ and set $\bx_1^n=(\bx_1,\dots,\bx_n)$.
A finite collection $f_1, \dots, f_N:\Rd \rightarrow \R$
  is called an {\em $\varepsilon$-- cover of $\F$ on $\bx_1^n$}
  if for any $f \in \F$ there exists  $i \in \{1, \dots, N\}$
  such that
  \[
\frac{1}{n} \sum_{k=1}^n |f(\bx_k)-f_i(\bx_k)| < \varepsilon.
  \]
  The {\em $\varepsilon$--covering number of $\F$ on $\bx_1^n$}
  is the  size $N$ of the smallest $\varepsilon$--cover
  of $\F$ on $\bx_1^n$ and is denoted by $\Nu_1(\varepsilon,\F,\bx_1^n)$.
  For $z \in \R$ and $\beta>0$ we define
  $T_\beta z = \max\{ - \beta,  \min\{z, \beta\}\}$.

\subsection{Outline of the paper}
The main result is presented and proven
in Section~\ref{se3}. In Section~\ref{se4} we prove a result
concerning the approximation of a smooth function on a manifold
by deep neural networks.

\section{Main result}
\label{se3}

Our main result is the following theorem, which presents a 
generalization 
bound of a least squares estimate based on fully connected DNNs in case 
when $\bX$ is concentrated on a $d^*$-dimensional Lipschitz-manifold.
\begin{theorem}
\label{th1}
Let $(\bX,Y), (\bX_1, Y_1), \dots, (\bX_n, Y_n)$
be independent and 
identically distributed random variables with values in $\Rd \times 
\R$ such that 
  \begin{equation*}
  \E\left\{ \exp(c_1 \cdot Y^2) \right\} < \infty
  \end{equation*}
  for some constant $c_1 > 0$.
  Let $p=q+s$ for some $q \in \N_0$ and $s \in (0,1]$,
  let $C>0$ and assume that the corresponding regression function
  $m(\cdot) = \EXP\{Y | X= \cdot\}$ is
  $(p,C)$-smooth and satisfies
  \[
\|m\|_{C^q(\Rd)} < \infty,
  \]
  and that the distribution of $\bX$ is
concentrated on a $d^*$-dimensional
  Lipschitz-manifold $\M$.
 Let $\tilde{m}_n$ be the least squares estimate defined by 
  \begin{equation*}
  \tilde{m}_n (\cdot) = \argmin_{h \in \mathcal{F}(L_{n},r_{n})} 
\frac{1}{n} \sum_{i=1}^n |Y_i - h(\bX_i)|^2
  \end{equation*}
  for some $L_n, r_n \in \N$,
  and define $m_n = T_{c_2 \cdot \log(n)} \tilde{m}_n$ for some $c_2 >0$ sufficiently large.

  Choose $c_{3}, c_{4} >0$ sufficiently large and set
  \[
  L_n = \left\lceil
    c_{3} \cdot \log n
    \right\rceil
    \quad \text{and} \quad 
    r_n = \left\lceil c_{4} \cdot
    n^{\frac{d^*}{2(2p+d^*)}} \right\rceil.
\]
 Then 
  \begin{equation*}
    \EXP \int |m_n(\bx) - m(\bx)|^2 {\PROB}_{\bX}(d\bx) \leq c_5 \cdot
    (\log n)^6
    \cdot
    n^{-\frac{2p}{2p+d^*}}
  \end{equation*}  
  holds for some constant $c_5 > 0$.
\end{theorem}

      \noindent
        {\bf Remark 1.} 
Let $\bv \in \R^{d-d^*}$.
Since
\[
[0,1]^{d^*} \times \{ \bv \}
\]
  is a $d^*$-dimensional Lipschitz manifold, it is easy to see
that Stone (1982) implies that the rate of convergence in 
Theorem~\ref{th1}
is optimal up to some logarithmic factor.\\

\noindent
{\bf Remark 2.}
The parameters $L_n$ and $r_n$ of the estimate in Theorem~\ref{th1}
depend on $d^*$ and $p$, which are usually unknown in practice.
But if they are chosen, e.g., by splitting of the sample
(cf., e.g., Chapter~7 in Gy\"orfi et al.\ (2002)), then the
corresponding
estimate, which does neither depend on $d^*$ or $p$, achieves
the same rate of convergence as our estimate in Theorem~\ref{th1}.\\

\noindent
{\bf Remark 3.}
       We conjecture that Theorem~\ref{th1} can also be extended to
        very deep fully connected neural network classes, where the
        number of neurons per layer is fixed and the number of hidden
        layers tends to infinity for sample size tending to infinity. 
In order to show this one could try to modify
the proof in Theorem~1b) in Kohler and Langer (2020). 
\\
\\
\noindent
{\bf Proof of Theorem \ref{th1}.}
Theorem 1 in Bagirov, Clausen and Kohler (2009)
(cf., Lemma 18 in Supplement B of Kohler and Langer (2020))
together with Lemma 19 in Supplement B of Kohler and Langer (2020) helps us to bound the expected $L_2$ error by 
\begin{align*}
&\mathbf{E} \int \left|m_n( \bx) - m( \bx)\right|^2 \PROB_\bX(d \bx)\\
& \leq \frac{c_{28} \cdot (\log(n))^2 c_{27} \cdot \log(n) \cdot \log(L_n \cdot r_n^2 ) \cdot L_n^2 \cdot r_n^2}{n}+ 2 \cdot \inf_{f \in \mathcal{F}(L_n,r_n)}
  \int\left|f(\bx) - m(\bx)\right|^2 \PROB_\bX(d \bx).
\end{align*}
Using this together with Theorem 2 below, where we choose
\begin{align*}
M = \lceil c_{33} \cdot n^{\frac{1}{2(2p+d^*)}} \rceil,
\end{align*}
shows the assertion.
\hfill $\Box$

\section{Approximation of smooth functions on a Lipschitz-manifold by 
deep neural networks}
\label{se4}

\subsection{An approximation result}
In this section we evaluate how well a DNN approximates a 
$(p,C)$-smooth 
function $f$ on a $d^*$-dimensional Lipschitz-manifold.
\begin{theorem}
  \label{th2} 
  Let $d \in \N$, let $d^* \in \{1, \dots, d\}$, let
  $\M$ be a $d^*$-dimensional Lipschitz-manifold, and let
  $1 \leq a < \infty$ such that $\M \subseteq [-a,a]^d$.
  Let $f:\Rd \rightarrow \R$ be $(p,C)$--smooth for some $p=q+s$,
  $q \in \N_0$, $s \in (0,1]$, and $C>0$.
    Let 
    $M \in \N$ be such that
    \begin{align*}
    M \geq 2 \ \text{and} 
    \ M^{2p} \geq c_6 \cdot \left(\max\left\{a, 
\|f\|_{C^q(\R^d)}\right\}\right)^{4(q+1)}
    \end{align*}
    holds for some sufficiently large constant $c_6 \geq 1$.
    Let $\sigma: \R \to \R$ be the ReLU activation function 
      \[
\sigma(x)= \max\{x,0\}.
\] 
There exists a neural network
\begin{align*}
\hat{f}_{\mathrm{net}} \in \mathcal{F}(L,r)
\end{align*} 
with
\begin{align*}
L= \lceil c_7 \cdot  \log(M) \rceil  \ \mbox{and} \   r = \lceil c_8 \cdot  M^{d^*} \rceil,
\end{align*}
such that
\begin{eqnarray}
 \| f-\hat{f}_\mathrm{net}\|_{\infty, \M} \leq
  c_{9} \cdot \left(\max\left\{1,  \|f\|_{C^q(\Rd)}
       \right\}\right)^{4(q+1)} \cdot M^{-2p}. 
  \label{th2eq1}
\end{eqnarray}
	\end{theorem}

\subsection{Idea of the proof of Theorem~\ref{th2}}
The proof of Theorem~\ref{th2} builds on the proof of Theorem~2a) in Kohler and 
Langer (2020). In particular,  we approximate a $(p,C)$-smooth function $f$ by a 
two-scale approximation, where we approximate piecewise Taylor 
polynomials with respect to a partition of
$\Rd$ into cubes with sidelength $1/M^{2}$. More precisely, we 
construct a coarse and fine grid of cubes with side length $1/M$ and 
$1/M^2$, respectively,
and approximate for the cube of the coarse grid,
which contains the $\bx$--value, 
the local Taylor polynomials on each cube
of the fine grid which is contained in the cube on
the coarse grid and which has non-empty intersection with $\M$.
\\
The difficulty compared to the result in Kohler and Langer (2020) is, that we only 
consider those cubes on our grid that have a non-empty intersection with our manifold.  This, in turn, means that our cubes are not necessarily next to each other.  It is therefore not possible to only use the information about a lower left corner of one cube to reproduce the whole grid as it was done in Kohler and Langer (2020).  Although our proof works somewhat differently, some parts are similar to the one of Theorem 2a) in Kohler and Langer (2020).  For convenience of the reader we will nevertheless present a complete proof, which sometimes means that we have to repeat arguments of Kohler and Langer (2020).

The partitions of $\Rd$ into the half-open equivolume cubes are defined as follows: Let
 \[
  C_{(k_1, \dots, k_d)}
  =
  \left[k_1 \cdot \frac{1}{M}, (k_1+1) \cdot \frac{1}{M}\right) \times \dots \times
    \left[k_d \cdot \frac{1}{M}, (k_d+1) \cdot \frac{1}{M}\right), \quad k_1, \dots, k_d \in \Z,
    \]
    and
    \[
    D_{(k_1, \dots, k_d)}=
    \left[k_1 \cdot \frac{1}{M^2}, (k_1+1) \cdot \frac{1}{M^2}\right) \times \dots \times
      \left[k_d \cdot \frac{1}{M^2}, (k_d+1) \cdot \frac{1}{M^2}\right), \quad k_1, \dots, k_d \in \Z,
  \]
  be the equivolume cubes with sidelengths $1/M$ and $1/M^2$, 
respectively. 
  Then the corresponding partitions are defined by 
  \begin{align}
\label{partition}
\mathcal{P}_1=\{C_{\bk} \, : \, \bk \in \Z^d \}
\quad \text{and} \quad
\mathcal{P}_2=\{D_{\bk} \, : \, \bk \in \Z^d\}.
\end{align}

%
For a partition of cubes $\P$ on $\Rd$, we denote by $C_{\P}(\bx)$ 
the cube $C$ that contains $\bx \in \Rd$.  The "bottom left" corner of 
some cube $C$ is denoted by $C_{\mathrm{left}}$. Therefore, one can 
describe 
the cube $C$ with side length $s$ (which is half-open as the cubes in 
$\P_1$ and $\P_2$) by a polytope as 
\begin{align*}
-x^{(j)} + C_{\mathrm{left}}^{(j)} \leq 0 \ 
\text{and} \ 
x^{(j)} - C_{\mathrm{\mathrm{left}}}^{(j)}-s < 0 \quad (j \in \{1, 
\dots, d\}).
\end{align*}
Furthermore, we describe by $C_{\delta}^0 \subset C$ the cube  
that contains all $\bx \in C$ that lie with a distance of at least 
$\delta$ to the borders of $C$, i.e., a polytope defined by
\begin{align*}
-x^{(j)} + C_{\mathrm{left}}^{(j)} \leq - \delta 
\ \text{and} \ 
x^{(j)} - C_{\mathrm{\mathrm{left}}}^{(j)}-s < -\delta \quad (j \in 
\{1, \dots, d\}).
\end{align*}
As every cube $C_{\bi} \in \mathcal{P}_1$ contains $M^d$ smaller 
cubes of $\mathcal{P}_2$, we denote those smaller cubes by 
$\tilde{C}_{1, \bi}, \dots, \tilde{C}_{M^d, \bi}$.
Here we order the cubes such that $\tilde{C}_{1, \bi}, \dots, \tilde{C}_{N_{\bi}, \bi}$ are all those
cubes which have a nonempty intersection with $\M$ (where
$N_\bi \in \{0,1, \dots, M^d\}$). We define by
\begin{align*}
T_{f,q,\bx_0}(\bx) = \sum_{\bj \in \N_0: \|\bj\|_1 \leq q} (\partial^{\bj} f)(\bx_0) \cdot \frac{(\bx-\bx_0)^{\bj}}{\bj!}
\end{align*}
the Taylor polynomial of total degree q around $\bx_0$.
\\
\\
Lemma~1 in Kohler (2014) implies that the piecewise Taylor polynomial
\begin{equation}
  \label{se4eq1}
T_{f,q,(C_{\P_2}(\bx))_{\mathrm{\mathrm{left}}}}(\bx)=\sum_{k \in \{1, 
\dots, M^d\}, \bi \in \Z^d } 
T_{f,q,(\tilde{C}_{k,\bi})_{\mathrm{\mathrm{left}}}}(\bx)
\cdot \mathds{1}_{\tilde{C}_{k,\bi}}(\bx) 
\end{equation}
satisfies
\[
  \Vert
f - \sum_{k \in \{1, \dots, M^d\}, \bi \in \Z^d } 
T_{f,q,(\tilde{C}_{k,\bi})_{\mathrm{\mathrm{left}}}}
\cdot \mathds{1}_{\tilde{C}_{k,\bi}}
\Vert_{\infty}
\leq
c_{10} \cdot (2 \cdot a \cdot d)^{2p} \cdot C \cdot \frac{1}{M^{2p}}.
\]
For $\bx \in \M$ we have
$\mathds{1}_{\tilde{C}_{k,\bi}}(\bx)=0$ if $C_\bi \cap \M = \emptyset$ or $k>N_\bi$,
hence
\[
T_{f,q,(C_{\P_2}(\bx))_{\mathrm{\mathrm{left}}}}(\bx)
=
\sum_{k \in \{1, \dots, N_\bi\}, \bi \in \Z^d : C_\bi \cap \M \neq 
\emptyset } T_{f,q,(\tilde{C}_{k,\bi})_{\mathrm{\mathrm{left}}}}(\bx)
\cdot \mathds{1}_{\tilde{C}_{k,\bi}}(\bx) 
\]
satisfies
\begin{eqnarray}
  &&
  \|
f - T_{f,q,(C_{\P_2}(\bx))_{\mathrm{\mathrm{left}}}}
\|_{\infty, \M}
=
  \|
  f -
  \sum_{k \in \{1, \dots, M^d\}, \bi \in \Z^d } 
T_{f,q,(\tilde{C}_{k,\bi})_{\mathrm{\mathrm{left}}}}(\bx)
\cdot \mathds{1}_{\tilde{C}_{k,\bi}}
\|_{\infty, \M}
\nonumber \\
&&
\leq
c_{10} \cdot (2 \cdot a \cdot d)^{2p} \cdot C \cdot \frac{1}{M^{2p}}.
\label{se3eq*}
\end{eqnarray}

If we use (\ref{se4eq1}) to approximate $f$ on some fixed compact set,
then it is easy to see that all summands except
some constant times
$M^{2d}$ of the summands
in (\ref{se4eq1}) are zero for all the $\bx$-values in the compact
set.
As our next lemma shows,
due to the fact that we use
$T_{f,q,(\tilde{C}_{k,\bi})_{\mathrm{\mathrm{left}}}}(\bx)$
only to approximate $f$ on our Lipschitz-manifold $\M$,
the number of summands in the definition
of
$T_{f,q,(\tilde{C}_{k,\bi})_{\mathrm{\mathrm{left}}}}(\bx)$, i.e.,
\[
\sum_{\bi \in \Z^d \, : \, C_\bi \cap \M \neq \emptyset}
N_\bi
=
\left|
\left\{
C \in \P_2 \, : \, C \cap \M \neq \emptyset
\right\}
\right|,
\]
is bounded by some constant times $M^{2d^*}$. Furthermore, we show
that $N_\bi \leq c_{11} \cdot M^{d^*}$ holds for all $\bi \in \Z^d$.

\begin{lemma}
  \label{le0}
  Let $\M$ be a $d^*$-dimensional Lipschitz-manifold.\\
 \noindent
  {\bf a)} Let $h \in (0,1]$ and set
      \[
      \P = \left\{
         [k_1 \cdot h, (k_1+1) \cdot h)
           \times \dots \times
                  [k_d \cdot h, (k_d+1) \cdot h)
                    \, : \,
                    k_1, \dots, k_d \in \Z
      \right\}.
      \]
      Then
      \[
      | \{ C \in \P \, : \, C \cap \M \neq \emptyset \} |
      \leq
      c_{12} \cdot
      \left( \frac{1}{h} \right)^{d^*},
      \]
      where
      $c_{12} = r \cdot (4 \cdot C_{\psi,2} \cdot \sqrt{d^*} +4)^{d^*}$.\\
\noindent      
      {\bf b)} Define $\P_1$, $\P_2$ and $N_\bi$ as above.  Then
\[
N_\bi \leq c_{13} \cdot M^{d^*}
\]
holds for all $\bi \in \Z^d$, where
$c_{13} = \max \{
1/C_{\psi,1}^{d^*},
3^{d^*} \cdot r^2 \cdot
(2 \cdot C_{\psi,2} \cdot \sqrt{d^*} +2)^{d^*}
\}
$.
  \end{lemma}
    {\bf Proof.} {\bf a)}
    Because of (\ref{de2eq2}) we have
    \begin{eqnarray*}
      &&
\M    \subseteq
    \bigcup_{j=1}^r
    \bigcup_{
k_1, \dots, k_{d^*} \in \{0,1, \dots, \lfloor 1/h \rfloor \}
    }
    \psi_j \left(
 [k_1 \cdot h, (k_1+1) \cdot h)
           \times \dots \times
                  [k_{d^*} \cdot h, (k_{d^*}+1) \cdot h)
    \right),
    \end{eqnarray*}
 hence we can bound 
\begin{align*}
&  | \{ C \in \P \, : \, C \cap \M \neq \emptyset \} |\\
&   \leq \sum_{j=1}^{r} \sum_{k_1=0}^{\lfloor \frac{1}{h}\rfloor} \dots \sum_{k_d^*=0}^{\lfloor \frac{1}{h}\rfloor} \\
& \quad |    \{ C \in \P \, : \, C \cap
\psi_j \left(
 [k_1 \cdot h, (k_1+1) \cdot h)
           \times \dots \times
                  [k_{d^*} \cdot h, (k_{d^*}+1) \cdot h)
                    \right)
                    \neq \emptyset \}|.
\end{align*}  
Consequently, it suffices to show that 
    \begin{eqnarray}
      \label{ple0eq1}
      &&
|    \{ C \in \P \, : \, C \cap
\psi_j \left(
 [k_1 \cdot h, (k_1+1) \cdot h)
           \times \dots \times
                  [k_{d^*} \cdot h, (k_{d^*}+1) \cdot h)
                    \right)
                    \neq \emptyset \}| \nonumber \\
                    &&
                    \leq (2 \cdot C_{\psi,2} \cdot \sqrt{d^*} +2)^{d^*}.
    \end{eqnarray}    
    The Lipschitz continuity of $\psi_j$ implies that
    \[
\psi_j \left(
 [k_1 \cdot h, (k_1+1) \cdot h)
           \times \dots \times
                  [k_{d^*} \cdot h, (k_{d^*}+1) \cdot h)
                    \right)
                    \]
                    is contained in a cube with sidelength
                    $2 \cdot C_{\psi,2} \cdot \sqrt{d^*} \cdot h$.
                    But any such cube has
                    a nonempty intersection with at most
                    \[
                    \left(
\frac{2 \cdot C_{\psi,2} \cdot \sqrt{d^*} \cdot h}{h} +2
\right)^{d^*}
= (2 \cdot C_{\psi,2} \cdot \sqrt{d^*}+2)^{d^*}
\]
many cubes from the partition $\P$. This shows the assertion.

\noindent
    {\bf b)} W.l.o.g. we can assume that $M \geq 1/C_{\psi,1}$.
    We have
    \begin{eqnarray*}
      N_\bi
      &=&
      | \{ 1 \leq j \leq M^d \, : \, \tilde{C}_{j,\bi} \cap \M \neq 0 \}|
      \\
      &=&
      \Bigg|
      \bigcup_{j=1}^{M^d}
      \bigcup_{l=1}^r
      \bigcup_{k_1=1}^{M-1}
      \dots
      \bigcup_{k_{d^*}=1}^{M-1}
      \\
      &&
      \quad
      \left\{ \tilde{C}_{j,\bi} \, : \, \tilde{C}_{j,\bi} \cap
\psi_l \left(
 \left[\frac{k_1}{M}, \frac{k_1+1}{M}\right)
           \times \dots \times
                  \left[\frac{k_{d^*}}{M} , \frac{k_{d^*}+1}{M}\right)
                    \right)
        \neq \emptyset \right\}\Bigg|.
      \end{eqnarray*}
Condition (\ref{de2eq1}) implies
    for any $\bx_1, \bx_2 \in [0,1]^{d^*}$ 
    \[
    \|\psi_l (\bx_1)- \psi_l (\bx_2)\| \geq  C_{\psi,1} \cdot
    \|\bx_1- \bx_2\| \geq \frac{1}{M} \cdot \|\bx_1- \bx_2\|.
    \]
    Using that two points in $ \tilde{C}_{j,\bi}$ have a supremum
    norm distance of at most $1/M^2$ this implies that
    for fixed $j \in \{1, \dots, M^d\}$ and
    $l \in \{1, \dots, r\}$
    there are at most $3^{d^*}$ different
    $(k_1, \dots, k_{d^*}) \in \{0,1, \dots, M-1\}^{d^*}$
    which satisfy
    \[
\tilde{C}_{j,\bi} \cap
\psi_l \left(
 \left[\frac{k_1}{M}, \frac{k_1+1}{M}\right)
           \times \dots \times
                  \left[\frac{k_{d^*}}{M} , \frac{k_{d^*}+1}{M}\right)
                    \right)
        \neq \emptyset.
    \]
    Using this we see
    \begin{eqnarray*}
    N_\bi &\leq& 3^{d^*} \cdot r \cdot
    \max_{l \in \{1, \dots, r\}, \atop
      k_1, \dots, k_{d^*} \in \{0, \dots, M-1\}}
    \\
    &&
    \quad
            \left|\left\{ 1 \leq j \leq M^d \, : \, \tilde{C}_{j,\bi} \cap
    \psi_j \left(
 \left[\frac{k_1}{M}, \frac{k_1+1}{M}\right)
           \times \dots \times
                  \left[\frac{k_{d^*}}{M} , \frac{k_{d^*}+1}{M}\right)
                    \right) \neq \emptyset
                    \right\} \right|.
    \end{eqnarray*}
    Using
    \begin{eqnarray*}
      &&
      \psi_j \left(
 \left[\frac{k_1}{M}, \frac{k_1+1}{M}\right)
           \times \dots \times
                  \left[\frac{k_{d^*}}{M} , \frac{k_{d^*}+1}{M}\right)
                    \right) \\
                    &&
                    \subseteq \bigcup_{j_1, \dots, j_{d^*} \in \{0, \dots, M-1\}} \\
                    &&
                      \psi_j \left(
 \left[\frac{k_1+j_1/M}{M}, \frac{k_1+(j_1+1)/M}{M}\right)
           \times \dots \times
                  \left[\frac{k_{d^*}+j_{d^*}/M}{M} ,
                    \frac{k_{d^*}+(j_{d^*}+1)/M}{M}\right)
                    \right)
      \end{eqnarray*}
    the assertion follows from (\ref{ple0eq1}).
    \hfill $\Box$
        
\noindent        
\\
To approximate $f(\bold{x})$ on $\M$ by neural networks our proof follows similar to the proof of Theorem 2a) in Kohler and Langer (2020) \textit{four} key steps:
\begin{enumerate}
\item[1.] Compute 
$T_{f,q,(C_{\P_2}(\bold{x}))_{\mathrm{left}}}(\bold{x})$ by using 
recursively defined functions.
\item[2.] Approximate the recursive functions by neural networks. The resulting network will be a good approximation for $f(\bold{x})$ in case that
\[\bold{x} \in \left( \bigcup_{\bk \in \Z^d } (D_{\bk})_{1/M^{2p+2}}^0 \right) \cap \M.
\]
\item[3.] Construct a neural network to approximate $w_{\P_2}(\bold{x}) \cdot f(\bold{x})$ for $\bx \in \M$, where
\begin{equation*}
w_{\P_2}(\bold{x}) = \prod_{j=1}^d \left(1- 2 \cdot M^2 \cdot 
\left|(C_{\mathcal{P}_{2}}(\bold{x}))_{\mathrm{left}}^{(j)} + 
\frac{1}{2 \cdot M^2} - x^{(j)}\right|\right)_+
\end{equation*}
is a linear tensorproduct B-spline
which takes its maximum value at the center of $C_{\P_{2}}(\bold{x})$, which
is nonzero in the inner part of $C_{\P_{2}}(\bold{x})$ and which
vanishes
outside of $C_{\P_{2}}(\bold{x})$. 
\item[4.] Apply those networks to $2^d$ slightly shifted partitions of $\P_2$ to approximate $f(\bold{x})$ in supremum norm.
\end{enumerate}

\subsection{Key step 1 of the proof of Theorem~\ref{th2}: A recursive 
definition of
  $T_{f,q,(C_{\P_2}(\bold{x}))_{\mathrm{left}}}(\bold{x})$ }
In the first key step we describe how to compute $T_{f,q,(C_{\P_2}(\bold{x}))_{\mathrm{left}}}(\bold{x})$ by recursively defined functions.  Those functions will later be approximated by neural networks.  \\
Assume $\bx \in \M$, and let  $\bi \in \Z^d$ such that we have $C_{\P_1}(\bold{x}) = C_{\bi}$. The recursion follows \textit{two} steps.
In a first step we compute the value of
$(\tilde{C}_{j,\bi})_{\mathrm{left}}$ for $j \in \{1, \dots, N_\bi\}$, the values of
$(\partial^{\bold{l}}f) ((\tilde{C}_{j,\bi})_{\mathrm{left}})$
for $j \in \{1, \dots, N_\bi\}$ and $\bl \in \N_0^d$ with $\|\bl\|_1 \leq q$ and  the length of the corresponding cubes $\tilde{C}_{j,\bi}$ for $j \in \{1, \dots, N_{\bi}\}$.
This can be done by computing the indicator function
$\mathds{1}_{C_{\bk}}$ multiplied by
$(\tilde{C}_{j,\bk})_{\mathrm{left}}$,
$(\partial^{\bold{l}} f)((\tilde{C}_{j,\bk})_{\mathrm{left}})$ and $1/M^2$
for each $\bk \in \Z^d$ with $C_\bk \cap \M \neq \emptyset$,
respectively.
Furthermore, we need the value of the input $\bold{x}$ in the further
recursive definition, therefore we shift this value by applying the
identity function. 
We set
\begin{align*}
\bm{\mathbf{\phi}}_{1,1} = (\phi_{1,1}^{(1)}, \dots, \phi_{1,1}^{(d)}) = \bold{x}, 
\end{align*}
\begin{align*}
  \bm{\phi}_{2, 1}^{(j)} &
  =
  (\phi_{2, 1}^{(j,1)}, \dots,
\phi_{2, 1}^{(j,d)}
  )
  = \sum_{\bk \in \Z^d : C_\bk \cap \M \neq \emptyset, \, j \leq N_\bk}
  (\tilde{C}_{j,\bk})_{\mathrm{left}}
  \cdot \mathds{1}_{C_{\bk}}(\mathbf{x}) 
\end{align*}
for $j \in \{1, \dots, \lceil c_{13} \cdot M^{d^*} \rceil\}$, 
\begin{align*}
  \phi_{3, 1}^{(\bl,j)} &= \sum_{\bk \in \Z^d : C_\bk \cap \M \neq \emptyset, \, j \leq N_\bk}
  (\partial^{\bold{l}} 
f)\left((\tilde{C}_{j,\bk})_{\mathrm{left}}\right)
  \cdot \mathds{1}_{C_{\bk}}(\mathbf{x}) 
\end{align*}
for $j \in \{1, \dots, \lceil c_{13} \cdot M^{d^*} \rceil \}$ and  $\mathbf{l} \in \N_0^d$ with $\|\bold{l}\|_1\leq q$, 
and 
\begin{align*}
\phi_{4,1}^{(j)} = \sum_{\bk \in \Z^d: C_{\bk} \cap \mathcal{M} \neq 0, j \leq N_{\bk}} \frac{1}{M^2} \cdot \mathds{1}_{C_{\bk}}(\mathbf{x}) 
\end{align*}
for $j \in \{1, \dots, \lceil c_{13} \cdot M^{d^*} \rceil\}$.
Here we have
 $\bm{\phi}_{2, 1}^{(j)}= \phi_{3,1}^{(\bl, j)}= \phi_{4, 1}^{(j)}=0$ for $j > N_\bi$.
\\
\\
Let $\bi \in \Z^d$
and
$j \in \{1, \dots, N_\bi\}$ such that $C_{\P_2}(\bold{x})=\tilde{C}_{j,\bi}$.
In a second step of the recursion we compute the value of  
$(C_{\P_2}(\bold{x}))_{\mathrm{left}}=(\tilde{C}_{j,\bi})_{\mathrm{
left}}$ and the values of $(\partial^{\bold{l}} 
f)\left((C_{\P_2}(\bold{x}))_{\mathrm{left}}\right)$ for $\bl \in 
\N_0^d$ with $\|\bl\|_1 \leq q$. It is easy to see
(cf.\ proof of Lemma~\ref{le1} below)
that each cube $\tilde{C}_{j,\bi}$
with $j \leq N_\bi$
can be defined by 
\begin{align}
\label{Aj}
\mathcal{A}^{(j)} = &\left\{\mathbf{x} \in \Rd: -x^{(k)} + \phi_{2,1}^{(j,k)}  \leq 0 \right. \notag\\
 & \hspace*{1.8cm} \left. \text{and} \ 
 x^{(k)} - \phi_{2,1}^{(j,k)}  
- \phi_{4,1}^{(j)} < 0 \ \text{for all} \ k \in \{1, \dots, d\}\right\}.
\end{align}
Thus, in our recursion we compute for each $j \in \{1, \dots, N_\bi\}$
the indicator function $\mathds{1}_{\mathcal{A}^{(j)}}$ multiplied by
$\bm{\mathbf{\phi}}_{2,1}^{(j)}$ or $\phi_{3, 1}^{(\bl, j)}$
for $\bl \in \N_0^d$ with $\|\bl\|_1 \leq q$.
Again we shift the value of $\bold{x}$ by applying the identity function.
We set
\begin{align*}
\bm{\mathbf{\phi}}_{1,2} = (\phi_{1,2}^{(1)}, \dots, \phi_{1,2}^{(d)})= \bm{\phi}_{1,1},
\end{align*}
\begin{align*}
  \bm{\phi}_{2,2}= (\phi_{2,2}^{(1)}, \dots, \phi_{2,2}^{(d)})=\sum_{j=1}^{
    \lceil c_{13} \cdot M^{d^*} \rceil
  } \bm{\phi}_{2,1}^{(j)}  \cdot \mathds{1}_{\mathcal{A}^{(j)}} \left(\bm{\phi}_{1,1}\right)
\end{align*}
and
\begin{align*}
\phi_{3,2}^{(\bl)} = \sum_{j=1}^{\lceil c_{13} \cdot M^{d^*} \rceil} \phi_{3, 1}^{(\bl, j)} \cdot \mathds{1}_{\mathcal{A}^{(j)}} \left(\bm{\phi}_{1,1}\right)
\end{align*}
for $\bl \in \N_0^d$ with $\|\bl\|_1 \leq q$.
In a last step, we compute the Taylor polynomial by
\begin{align*}
\phi_{1,3} = &\sum_{\substack{\bj \in \N_0: \|\bj\|_1 \leq q}} \frac{\phi_{3, 2}^{(\bj)}}{\bj!} \cdot \left(\bm{\phi}_{1,2} - \bm{\phi}_{2,2}\right)^{\bj}.
\end{align*}
Our next lemma shows that this recursion computes our piecewise Taylor polynomial on $\M$.

\begin{lemma}
\label{le1}
Let $p=q+s$ for some $q \in \N_0$ and $s \in (0,1]$, and let $C > 0$.
Let $f: \Rd \to \R$ be a $(p,C)$-smooth function and let 
$T_{f,q,(C_{\mathcal{P}_2}(\bold{x}))_{\mathrm{left}}}$ be the Taylor 
polynomial of total degree $q$ around 
$(C_{\mathcal{P}_2}(\bold{x}))_{\mathrm{left}}$. Define $\phi_{1,3}$ 
recursively as above. Then we have for any $\bx \in \M$:
  \[
\phi_{1,3}=T_{f,q,(C_{\mathcal{P}_2}(\bold{x}))_{\mathrm{left}}}(\bold
{x}).
  \]
\end{lemma}

\noindent
{\bf Proof.}
Let $\bx \in \M$ and let $\bi \in \Z^d$ and $j \in \{1, \dots, N_\bi\}$
be such that we have $C_{\P_2(\bx)}=\tilde{C}_{j,\bi}$. Then we have
$\bx \in C_\bi$, and our definitions above imply
\[
\bm{\phi}_{1,1}=\bx, \,
\quad
\bm{\phi}_{2,1}^{(k)}
=
(\tilde{C}_{k,\bi})_{\mathrm{left}} \cdot \mathds{1}_{\{k \leq 
N_\bi\}}(\bx),
\quad
\phi_{3,1}^{(\bl,k)}
=
(\partial^{\bold{l}} 
f)\left((\tilde{C}_{k,\bi})_{\mathrm{left}}\right)
\cdot \mathds{1}_{\{k \leq N_\bi\}}(\bx)
\]
and
\[
\phi_{4,1}^{(k)} = \frac{1}{M^2} \cdot \mathds{1}_{\{k \leq N_{\bi}\}}(\bx)
\]
for all $k \in \{1, \dots, \lceil c_{13} \cdot M^{d^*} \rceil 
\}$ and $\bl \in \N_0^d$ with $\|\bl\|_1 \leq q$.
This implies
\[
\A^{(k)}=\tilde{C}_{k,\bi} \quad \text{for } k \leq N_\bi,
\]
and using again our definitions above we see
\[
\bm{\phi}_{1,2}=\bx, \quad 
\bm{\phi}_{2,2}=
\sum_{k=1}^{N_\bi}
(\tilde{C}_{k,\bi})_{left}
\cdot
\mathds{1}_{ \tilde{C}_{k,\bi} }(\bx)
=
(\tilde{C}_{j,\bi})_{\mathrm{left}}
\]
and
\[
\phi_{3,2}^{(\bl)}
=
\sum_{k=1}^{N_\bi}
(\partial^{\bl} f)
(
(\tilde{C}_{k,\bi})_{left}
)
\cdot
\mathds{1}_{ \tilde{C}_{k,\bi} }(\bx)
=
(\partial^{\bold{l}} 
f)\left((\tilde{C}_{j,\bi})_{\mathrm{left}}\right).
\]
Consequently, it holds
\begin{eqnarray*}
\phi_{1,3} &=&
\sum_{\substack{\bj \in \N_0: \|\bj\|_1 \leq q}} \frac{\phi_{3, 2}^{(\bj)}}{\bj!} \cdot \left(\bm{\phi}_{1,2} - \bm{\phi}_{2,2}\right)^{\bj}
=
\sum_{\substack{\bj \in \N_0: \|\bj\|_1 \leq q}} \frac{
  (\partial^{\bold{l}} 
f)\left((\tilde{C}_{j,\bi})_{\mathrm{left}}\right)
}{\bj!}
\cdot \left(\bx -  (\tilde{C}_{j,\bi})_{\mathrm{left}}\right)^{\bj}\\
&=&
T_{f,q,(C_{\mathcal{P}_2}(\bold{x}))_{\mathrm{left}}}(\bold{x}).
\end{eqnarray*}
\hfill $\Box$

\subsection{Key step 2 of the proof of Theorem~\ref{th2}: 
Approximating $\phi_{1,3}$ by neural networks}
In key step 2 we approximate the functions $\bm{\phi}_{1,1}$, $\bm{\phi}_{2,1}^{(j)}$, $\phi_{3,1}^{(\bl, j)}$,  $\phi_{4,1}^{(j)}$, 
 $\bm{\phi}_{1,2}$, $\bm{\phi}_{2,2}$, $\phi_{3, 2}^{(\bl)}$, 
$\phi_{1,3}$ ($j \in \{1, \dots,
\lceil c_{13} \cdot M^{d^*} \rceil \}$, $\bl \in \N_0^d$ with $\|\bl\|_1 \leq q$) by neural networks.  
By using the following two computing operations for neural networks, we can combine smaller networks in one large neural network:
\\
\\
\textit{Combined neural network:} Let $f \in \mathcal{F}(L_f, r_f)$ 
and $g \in \mathcal{F}(L_g,r_g)$ with $L_f, L_g, r_f, r_g \in \N$, 
then we call $f \circ g$ the \textit{combined network}, which is 
contained in the network class $\mathcal{F}(L_f+L_g, \max\{r_f, r_g))$. 
Here, the output of the network $g$ is the input 
of the network $f$ and the total number of hidden layers equals the 
sum of the hidden layers of both networks $f$ and $g$ (cf.\ Figure~3 
in Kohler and Langer (2020)).
\\
\\
\textit{Parallelized neural network:} Let $f \in \mathcal{F}(L, r_f)$ 
and $g \in \mathcal{F}(L, r_g)$ 
be two networks with the same number of hidden layers $L \in \N$. Then 
we call $(f,g) \in \mathcal{F}(L, r_f+r_g)$ the parallelized network, 
which computes $f$ and $g$ in parallel in a joint network.  
\\
\\
The final network of this step approximates $f(\bx)$ in case that $\bold{x} \in \M$ does not lie close to the boundary of any cube of $\P_2$, i.e., for
\begin{align*}
  \bold{x} \in  \left(
  \bigcup_{\bk \in \Z^d}(D_{\bk})_{1/M^{2p+2}}^0
  \right)
  \cap \M .
\end{align*}

\begin{lemma}
\label{le2} 
Let $\sigma:\R \to \R$ be the ReLU activation function 
$\sigma(x) = \max\{x,0\}$. Let $\mathcal{P}_2$ be defined as in 
\eqref{partition}. Let $p = q+s$ for some $q \in \N_0$ and $s \in 
(0,1]$, and let $C >0$. Let $f: \Rd \to \R$ be a $(p,C)$-smooth 
function.
                 Let
                 $\M$ be a $d^*$-dimensional Lipschitz-manifold, and let
                 $1 \leq a < \infty$ such that
$\M \subseteq [-a,a]^d$.
Then there exists for $M \in \N$ with
    \begin{align*}
    M^{2p} \geq c_{14} \cdot \left(\max\left\{3a, \|f\|_{C^q(\R^d)}\right\}\right)^{4(q+1)}
    \end{align*}
 a neural network
$\hat{f}_{\P_2} \in \mathcal{F}(L,r)$ with
\begin{itemize}
\item[(i)] $L= 4+\lceil \log_4(M^{2p})\rceil \cdot \lceil \log_2(\max\{q+1,2 \})\rceil$
\item[(ii)] $r=\max\left\{\left(\binom{d+q}{d} + d\right) \cdot
  \lceil c_{13} \cdot M^{d^*} \rceil \cdot 2 \cdot (2+2d)+2d, 18 \cdot (q+1) \cdot \binom{d+q}{d}\right\}$
\end{itemize}
such that 
\begin{align*}
&|\hat{f}_{\mathcal{P}_2}(\bold{x}) - f(\bold{x})| \leq c_{15} \cdot \left(\max\left\{3a, \|f\|_{C^q(\R^d)}\right\}\right)^{4(q+1)} \cdot \frac{1}{M^{2p}}
\end{align*}
holds for all $\bold{x} \in
\left( \bigcup_{\bk \in \Z^d} \left(D_{\bk}\right)_{1/M^{2p+2}}^0
\right) \cap \M$.
The network value is bounded by 
\begin{align*}
|\hat{f}_{\mathcal{P}_2}(\bold{x})| \leq 2 \cdot e^{4ad} \cdot \max\left\{\|f\|_{C^q(\R^d)}, 1\right\}
\end{align*}
for all $\bold{x} \in \R^d$.
\end{lemma}

In the proof of Lemma~\ref{le2} we will need several
auxiliary neural networks, which we introduce next:\\
\\
\textit{Identity network:} As we are using the ReLU activation 
function, we can exploit its projection property to shift input values 
in the next hidden layer or to synchronize the number 
of hidden layers for two networks, which are computed in parallel.  
Here we use the network $\hat{f}_{\mathrm{id}}: \R \to \R$,
\begin{align*}
\hat{f}_{\mathrm{id}}(z) = \sigma(z) - \sigma(-z) = z, 
\quad z \in \R, 
\end{align*}
and 
\begin{align*}
\hat{f}_{\mathrm{id}}(\bold{x}) = 
\left(\hat{f}_{\mathrm{id}}\left(x^{(1)}\right), \dots, 
\hat{f}_{\mathrm{id}}\left(x^{(d)}\right)\right) = \left(x^{(1)}, 
\dots, x^{(d)}\right), \quad \bold{x} \in \Rd.
\end{align*}
Furthermore, we will use the abbreviations
\begin{align*}
&\hat{f}_{\mathrm{id}}^0(\bold{x}) = \bold{x}, \quad \bold{x} \in 
\R^d\\
&\hat{f}_{\mathrm{id}}^{t+1}(\bold{x}) = 
\hat{f}_{\mathrm{id}}\left(\hat{f}_{\mathrm{id}}^t(\bold{x})\right) = 
\bold{x}, \quad t \in \N_0, \bold{x} \in \R^d.
\end{align*}
\noindent
\textit{Network for polynomials:} Let $\mathcal{P}_N$ be the linear span of all monomials of the form 
\begin{align*}
\prod_{k=1}^d \left(x^{(k)}\right)^{r_k}
\end{align*}
for some $r_1, \dots, r_d \in \N_0$, $r_1+\dots+r_d \leq N$. Then, $\mathcal{P}_N$ is a linear vector space of functions of dimension 
\begin{align*}
\dim \ \mathcal{P}_N = \left|\left\{(r_0, \dots, r_d) \in \N_0^{d+1}: r_0+\dots+r_d = N \right\}\right| = \binom{d+N}{d}.
\end{align*}
The next lemma describes a neural network that approximates functions 
of the class $\mathcal{P}_N$ multiplied by an additional factor. This 
modified form of polynomials 
is later needed in the construction of our network of Lemma~\ref{le2}.
\begin{lemma}
\label{le4} Let $a \geq 1$.
Let $m_1, \dots, m_{\binom{d+N}{d}}$ denote all monomials in 
$\mathcal{P}_N$ for some $N \in \N$. Let $r_1, \dots, 
r_{\binom{d+N}{d}} \in \R$, define 
\begin{align*}
p\left(\bold{x}, y_1, \dots, y_{\binom{d+N}{d}}\right) = 
\sum_{i=1}^{\binom{d+N}{d}} r_i \cdot y_i \cdot m_i(\bold{x}), \quad 
\bold{x} \in [-a,a]^d, y_i \in [-a,a],
\end{align*}
and set $\bar{r}(p) = \max_{i \in \left\{1, \dots, \binom{d+N}{d}\right\}} |r_i|$. Let $\sigma: \mathbb{R} \to \R$ be the ReLU activation function $\sigma(x) = \max\{x,0\}$. Then for any  
\begin{equation}
\label{le3eq1}
R \geq \log_4 (2 \cdot 4^{2 \cdot (N+1)} \cdot a^{2 \cdot (N+1)})
\end{equation}
 a neural network 
\begin{align*}
\hat{f}_{p} \in \mathcal{F}(L,r)
\end{align*}
with $L=R \cdot \lceil\log_2(N+1)\rceil$ and $r =18 \cdot (N+1) \cdot \binom{d+N}{d}$ exists, such that
\begin{align*}
\left|\hat{f}_{p}\left(\bold{x}, y_1, \dots, y_{\binom{d+N}{d}}\right) - p\left(\bold{x}, y_1, \dots, y_{\binom{d+N}{d}}\right) \right| \leq c_{16} \cdot \bar{r}(p) \cdot a^{4(N+1)} \cdot 4^{-R}
\end{align*}
for all $\bold{x} \in [-a,a]^d$, $y_1, \dots, y_{\binom{d+N}{d}} \in 
[-a,a]$, 
where $c_{16}$ depends on $d$ and $N$.
\end{lemma}
\noindent
{\bf Proof.}
See Lemma~5 in Supplement~A of Kohler and Langer (2020).
\hfill $\Box$
\\
\\
\noindent
\textit{Network for multidimensional indicator functions:} 
The next lemma presents a network that approximates the 
multidimensional indicator function and the multidimensional indicator 
function multiplied by an additional factor. 
\begin{lemma}
\label{le5}
Let $\sigma: \R \to \R$ be the ReLU activation function $\sigma(x) = \max\{x,0\}$. Let $R \in \N$. Let $\mathbf{a}, \mathbf{b} \in \Rd$ with
\begin{align*}
b^{(i)} - a^{(i)} \geq \frac{2}{R} \ \text{for all} \ i \in \{1, 
\dots, d\}
\end{align*}
and let
\begin{align*}
&K_{1/R} = \big\{\bold{x} \in \Rd: x^{(i)} \notin [a^{(i)}, a^{(i)}+1/R) \cup (b^{(i)} - 1/R, b^{(i)})\\
& \hspace*{8cm}  \text{for all} \ i \in \{1, \dots, d\}\big\}.
\end{align*}
a) Then the network
\begin{align*}
\hat{f}_{\mathrm{ind}, [\bold{a}, \bold{b})}(\bold{x}) &= 
\sigma\bigg(1-R \cdot \sum_{i=1}^d \left(\sigma\left(a^{(i)} + 
\frac{1}{R} - x^{(i)}\right)\right.\\
& \hspace*{3cm} \left. + \sigma\left(x^{(i)} - b^{(i)} + \frac{1}{R}\right) \right)\bigg)
\end{align*}
of the class $\mathcal{F}(2, 2d)$ satisfies
\begin{align*}
\hat{f}_{\mathrm{ind}, [\bold{a}, \bold{b})}(\bold{x}) = \mathds{1}_{ 
[\bold{a}, \bold{b})}(\bold{x})
\end{align*}
for $\bold{x} \in K_{1/R}$ and 
\begin{align*}
\left|\hat{f}_{\mathrm{ind}, [\bold{a}, \bold{b})}(\bold{x}) - 
\mathds{1}_{[\bold{a}, \bold{b})}(\mathbf{x})\right| \leq 1
\end{align*}
for $\bold{x} \in \Rd$.
\\
b) Let $|s| \leq R$. Then the network
\begin{align*}
\hat{f}_{\mathrm{test}}(\bold{x}, \mathbf{a}, \mathbf{b}, s) &= 
\sigma\bigg(\hat{f}_{\mathrm{id}}(s)-R^2 \cdot \sum_{i=1}^d 
\left(\sigma\left(a^{(i)} + \frac{1}{R} - x^{(i)}\right)\right.\\
& \left. \hspace*{4cm} + \sigma\left(x^{(i)} - b^{(i)} + \frac{1}{R}\right)\right)\bigg)\\
& \quad - \sigma\bigg(-\hat{f}_{\mathrm{id}}(s)-R^2 \cdot \sum_{i=1}^d 
\left(\sigma\left(a^{(i)} + \frac{1}{R} - x^{(i)}\right)\right.\\
& \left. \hspace*{4cm} + \sigma\left(x^{(i)} - b^{(i)} + \frac{1}{R}\right)\right)\bigg)
\end{align*}
of the class $\mathcal{F}(2, 2 \cdot (2d+2))$
satisfies
\begin{align*}
  \hat{f}_{\mathrm{test}}(\bold{x}, \mathbf{a}, \mathbf{b}, s)
  =
  s \cdot \mathds{1}_{[\bold{a}, \bold{b})}(\bold{x})
\end{align*}
for $\bold{x} \in K_{1/R}$ and 
\begin{align*}
  \left|\hat{f}_{\mathrm{test}}(\bold{x}, \mathbf{a}, \mathbf{b}, s) -
  s \cdot \mathds{1}_{[\bold{a}, \bold{b})}(\bold{x})\right| \leq |s|
\end{align*}
for $\bold{x} \in \Rd$.
\end{lemma}
\noindent
{\bf Proof.} 
See Lemma~6  in Supplement~A of Kohler and Langer (2020).
\hfill $\Box$
\\
\\
{\bf Proof of Lemma \ref{le2}.}
In a \textit{first step of the proof} we describe how the recursively 
defined function $\phi_{1,3}$  
can be approximated by neural networks.  To approximate an indicator function
$\mathds{1}_{[\bold{a}, \bold{b})}(\bold{x})$ for some $\mathbf{a}, 
\mathbf{b} \in \Rd$ and $B_M \in \N$ with
\begin{align*}
b^{(i)} - a^{(i)} \geq \frac{2}{B_M} \ \text{for all} \ i \in \{1, 
\dots, d\}
\end{align*}
we will 
use the network
\begin{align*}
\hat{f}_{\mathrm{ind}, [\bold{a}, \bold{b})} \in \mathcal{F}(2, 2d)
\end{align*}
of Lemma~\ref{le5} (with $R=B_M$).  With the networks
\begin{align*}
\hat{f}_{\mathrm{test}} \in \mathcal{F}(2, 2 \cdot (2d+2))
\end{align*}
of Lemma~\ref{le5}
(again with $R=B_M$)
we approximate functions of the form
\begin{align*}
s \cdot \mathds{1}_{[\bold{a}, \bold{b})}(\bold{x}).
\end{align*}
Lemma \ref{le5} implies that for $|s| \leq B_M$ and
$\bx \in \Rd$ with
\begin{align*}
x^{(i)} \notin \Big[a^{(i)}, a^{(i)} + \frac{1}{B_M}\Big) \cup 
\Big(b^{(i)} - \frac{1}{B_M}, b^{(i)} \Big) \ \text{for all} \ i \in 
\{1, \dots, d\}
\end{align*}
we have
\begin{align*}
\hat{f}_{\mathrm{ind}, [\bold{a}, \bold{b})}(\bold{x}) = 
\mathds{1}_{[\bold{a}, \bold{b})}(\bold{x})
\end{align*}
and
\begin{align*}
\hat{f}_{\mathrm{test}}(\bold{x}, \mathbf{a}, \mathbf{b}, s)(\bold{x}) 
= s \cdot \mathds{1}_{[\bold{a}, \bold{b})}(\bold{x}).
\end{align*}
For some vector $\bold{v} \in \R^d$ it follows
\begin{align*}
&\bold{v} \cdot \hat{f}_{\mathrm{ind}, [\bold{a}, \bold{b})}(\bold{x}) 
= \left(v^{(1)} \cdot \hat{f}_{\mathrm{ind}, [\bold{a}, 
\bold{b})}(\bold{x}), \dots, v^{(d)} \cdot \hat{f}_{\mathrm{ind}, 
[\bold{a}, \bold{b})}(\bold{x})\right).
\end{align*}
To compute the final Taylor polynomial in $\phi_{1,3}$ we use the network
\begin{align*}
\hat{f}_{p} \in \mathcal{F}\left(B_{M,p} \cdot \lceil \log_2(\max\{q+1, 2\}) \rceil, 18 \cdot (q+1) \cdot \binom{d+q}{d}\right)
\end{align*} 
from Lemma~\ref{le4}
(with $R=B_{M,p}$)
satisfying 
\begin{align}
\label{fpeq}
&\left|\hat{f}_{p}\left(\bold{z}, y_1, \dots, y_{\binom{d+q}{q}}\right) - p\left(\bold{z}, y_1, \dots, y_{\binom{d+q}{q}}\right)\right| \notag\\
& \leq c_{16} \cdot \bar{r}(p)
\cdot \left(\max\left\{3a, \|f\|_{C^q(\R^d)}\right\}\right)^{4(q+1)} \cdot 4^{-B_{M,p}}
\end{align}
for all $z^{(1)}, \dots, z^{(d)}, y_1, \dots, y_{\binom{d+q}{d}}$ contained in
\begin{align*}
\left[-\max\left\{3a, \|f\|_{C^q(\R^d)}\right\}, \max\left\{3a, \|f\|_{C^q(\R^d)}\right\}\right],
\end{align*}
where $B_{M,p} \in \N$ satisfying
\begin{align*}
  B_{M,p} \geq \log_4 \left(
  \max \{ c_{16}, 2 \cdot 4^{2 \cdot (q+1)} \}
  \cdot \left(\max\left\{3a, \|f\|_{C^q(\R^d)}\right\}\right)^{2 \cdot (q+1)}\right)
\end{align*}
is properly chosen (cf.\ (\ref{le3eq1})).
In case that $q=0$ we use a polynomial of degree $1$ 
where the $r_i$'s of all coefficients greater than zero are chosen as 
zero. That is why we changed $\log_2(q+1)$ to $\log_2(\max\{q+1,2\})$ 
in the definition of $L$ in
Lemma~\ref{le4}.
\\
\\
Each network of the recursion of $\phi_{1,3}$ is now computed by a neural network. To compute the values of $\bm{\phi}_{1,1}$, $\bm{\phi}_{2,1}^{(j)}$,
$\phi_{3, 1}^{(\bl, j)}$
and $\phi_{4, 1}^{(j)}$ we use for $j \in \{1, \dots,\lceil c_{13} \cdot  M^{d^*} \rceil \}$ and $\mathbf{l} \in \N_0^d$ with $\|\bl\|_1 \leq q$ the networks
\[
\bm{\hat{\phi}}_{1,1} = \left(\hat{\phi}_{1,1}^{(1)}, \dots, 
\hat{\phi}_{1,1}^{(d)}\right) = \hat{f}_{\mathrm{id}}^2 (\bold{x}),
\]
\[
\bm{\hat{\phi}}_{2,1}^{(j)} =
(\hat{\phi}_{2,1}^{(j,1)}, \dots, \hat{\phi}_{2,1}^{(j,d)}) = \sum_{\bi \in \Z^d \, : \, C_\bi \cap \M \neq \emptyset,
  \, j \leq N_\bi}
(\tilde{C}_{j,\bi})_{\mathrm{left}} \cdot 
\hat{f}_{\mathrm{ind},{C_{\bi}}}(\bold{x}),
\]
\[
\hat{\phi}_{3,1}^{(\mathbf{l},j)} = \sum_{\bi \in \Z^d \, : \, C_\bi \cap \M \neq \emptyset,
  \, j \leq N_\bi}
(\partial^{\bl} f)\left((\tilde{C}_{j,i})_{\mathrm{left}}\right) \cdot 
\hat{f}_{\mathrm{ind},{C_{\bi}}}(\bold{x}).
\]
and
\[
\hat{\phi}_{4,1}^{(j)} =  \sum_{\bi \in \Z^d \, : \, C_\bi \cap \M \neq \emptyset,
  \, j \leq N_\bi} \frac{1}{M^2} \cdot 
\hat{f}_{\mathrm{ind},{C_{\bi}}}(\bold{x}).
\]
To compute $\bm{\phi}_{1,2}$, $\bm{\phi}_{2,2}$ and $\phi_{3, 2}^{(\bl)}$
we use the networks
\[
\bm{\hat{\phi}}_{1,2}= \left(\hat{\phi}_{1,2}^{(1)}, \dots, 
\hat{\phi}_{1,2}^{(d)}\right) = \hat{f}_{\mathrm{id}}^{2} 
(\bm{\hat{\phi}}_{1,1}),
\]
\begin{align}
\label{neur31}
\hat{\phi}_{2,2}^{(k)} = \sum_{j=1}^{\lceil c_{13} \cdot M^{d^*} 
\rceil} \hat{f}_{\mathrm{test}}\left(\bm{\hat{\phi}}_{1,1}, 
\bm{\hat{\phi}}_{2,1}^{(j)}, \bm{\hat{\phi}}_{2,1}^{(j)} 
+\hat{\phi}_{4,1}^{(j)} \cdot \mathbf{1}, \hat{\phi}_{2,1}^{(j,k)}\right)
\end{align}
for $k \in \{1, \dots, d\}$ and
\begin{align*}
\bm{\hat{\phi}}_{2,2}=(\hat{\phi}_{2,2}^{(1)}, \dots, \hat{\phi}_{2,2}^{(d)})
\end{align*}
and
\begin{align}
\label{neur32}
&\hat{\phi}_{3, 2}^{(\bl)} = \sum_{j=1}^{\lceil c_{13} \cdot M^{d^*} \rceil}
\hat{f}_{\mathrm{test}}\left(\bm{\hat{\phi}}_{1,1}, 
\bm{\hat{\phi}}_{2,1}^{(j)},
\bm{\hat{\phi}}_{2,1}^{(j)} +\hat{\phi}_{4,1}^{(j)} \cdot \mathbf{1}, \hat{\phi}_{3, 1}^{(\bl, j)}\right).
\end{align}
%
Choose $\bl_1, \dots, \bl_{\binom{d+q}{d}}$ such that
\begin{align*}
\left\{\bl_1, \dots, \bl_{\binom{d+q}{d}}\right\} = \left\{(s_1, \dots, s_d) \in \N_0^d: s_1+\dots+s_d \leq q \right\}
\end{align*}
holds. 
The value of $\phi_{1,3}$ can then be computed by 
\begin{align}
\label{fp}
\hat{\phi}_{1,3} = \hat{f}_p\left(\bold{z}, y_1, \dots, y_{\binom{d+q}{d}}\right),
\end{align}
where 
\begin{align*}
\bold{z}= \bm{\hat{\phi}}_{1,2} - \bm{\hat{\phi}}_{2,2}
\end{align*}
and 
\begin{align*}
y_v = \hat{\phi}_{3, 2}^{(\bl_v)} 
\end{align*}
for $v \in \left\{1, \dots, \binom{d+q}{d}\right\}$. 
The coefficients $r_1, \dots, r_{\binom{d+q}{d}}$ in Lemma~\ref{le4} 
are chosen as 
\begin{align*}
r_i = \frac{1}{\bl_i!}, \quad i \in \left\{1, \dots, \binom{d+q}{d}\right\}.
\end{align*}

It is easy to see that the network $\hat{\phi}_{1,3}$ forms a composed 
network, where the networks $\bm{\hat{\phi}}_{1,1}$,
$\bm{\hat{\phi}}_{2, 1}^{( 1)}$, \dots, $\bm{\hat{\phi}}_{2,1}^{(\lceil c_{13} \cdot M^{d^*} \rceil)}$,
$\hat{\phi}_{3, 1}^{(\bl_v, 1)}$, \dots,
$\hat{\phi}_{3,1}^{(\bl_v, \lceil c_{13} \cdot M^{d^*} \rceil)}$,$\hat{\phi}_{4, 1}^{(1)}$, \dots,
$\hat{\phi}_{4,1}^{(\lceil c_{13} \cdot M^{d^*} \rceil)}$
and the networks $\bm{\hat{\phi}}_{1,2}, 
\bm{\hat{\phi}}_{2,2}, \hat{\phi}_{3, 2}^{(\bl_v)}$ 
$(v \in \{1, \dots, \binom{d+q}{d}\})$ are computed in parallel (i.e., 
in the same layers),
respectively. Thus, we can conclude that
\begin{align*}
  (\bm{\hat{\phi}}_{1,1},   \bm{\hat{\phi}}_{2, 1}^{(1)}, \dots, \bm{\hat{\phi}}_{2, 1}^{(\lceil c_{13} \cdot M^{d^*} \rceil)},
  \hat{\phi}_{3, 1}^{(\bl_v,1)}, \dots, \hat{\phi}_{3, 1}^{(\bl_v, \lceil c_{13} \cdot M^{d^*} \rceil)}, \hat{\phi}_{4,1}^{(1)}, \dots, \hat{\phi}_{4,1}^{(\lceil c_{13} \cdot M^{d^*} \rceil)}
  ) 
\end{align*}
needs $L_1=2$ hidden layers and $r_1=2d+\left(d+ \binom{d+q}{d}+1\right) \cdot
\lceil c_{13} \cdot M^{d^*} \rceil
\cdot 2d$ 
neurons per layer in total.
\\
\\
Furthermore, the parallelized network
\begin{align*}
(\bm{\hat{\phi}}_{1,2}, \bm{\hat{\phi}}_{2,2}, \hat{\phi}_{3, 2}^{(\bl_v)}) 
\end{align*}
needs $L_2=L_1+2=4$ hidden layers and 
\begin{align*}
r_2&=\max\left\{r_1, 2d+d \cdot \lceil c_{13} \cdot M^{d^*} \rceil
\cdot 2 \cdot (2d+2)+\binom{d+q}{d} \cdot
\lceil c_{13} \cdot M^{d^*} \rceil \cdot 2 \cdot (2d+2)\right\}\\
&= 2d+\left(d +\binom{d+q}{d}\right)\cdot \lceil c_{13} \cdot M^{d^*} \rceil
 \cdot 2 \cdot (2d+2)
\end{align*}
neurons per layer.
Finally we have that $\hat{\phi}_{1,3}$ lies in the class
\begin{align*}
\mathcal{F}\left(4+B_{M,p} \cdot \lceil \log_2(\max\{q+1, 2\}) \rceil, r\right)
\end{align*}
with 
\begin{align*}
r=\max\left\{r_2, 18 \cdot (q+1) \cdot \binom{d+q}{d}\right\}.
\end{align*}
Here we have used that
\begin{align*}
\mathcal{F}(L, r') \subseteq \mathcal{F}(L, r)
\end{align*}
for $r' \leq r$.
We set
\begin{align*}
\hat{f}_{\mathcal{P}_2}(\bold{x}) = \hat{\phi}_{1,3}.
\end{align*}
In a \textit{second step of the proof} we analyze the error of the network $\hat{f}_{\mathcal{P}_2}$ in case that 
\begin{align*}
B_M \geq M^{2p+2} 
\end{align*}
and 
\begin{align*}
  \bold{x} \in
  \left(
  \bigcup_{\bk \in \Z^d} \left(D_{\bk}\right)_{1/M^{2p+2}}^0
  \right) \cap \M.
\end{align*}
From Lemma~\ref{le5} we can conclude that the networks 
$\bm{\hat{\phi}}_{1,1}$,
$\bm{\hat{\phi}}_{2, 1}^{( 1)}$, \dots, $\bm{\hat{\phi}}_{2,1}^{(\lceil c_{13} \cdot M^{d^*} \rceil)}$,
$\hat{\phi}_{3, 1}^{(\bl_v, 1)}$, \dots,
$\hat{\phi}_{3,1}^{(\bl_v, \lceil c_{13} \cdot M^{d^*} \rceil)}$, $\hat{\phi}_{4,1}^{(1)}, \dots, \hat{\phi}_{4,1}^{(\lceil c_{13} \cdot M^{d^*}\rceil)}$ and the networks $\bm{\hat{\phi}}_{1,2}, \bm{\hat{\phi}}_{2,2}, \hat{\phi}_{3, 2}^{(\bl_v)}$ $(v \in \{1, \dots, \binom{d+q}{d}\})$ compute the corresponding functions $\bm{\phi}_{1,1}$,
$\bm{\phi}_{2,1}^{(1)}$, \dots,
$\bm{\phi}_{2,1}^{(\lceil c_{13} \cdot M^{d^*} \rceil)}$,
$\phi_{3, 1}^{(\bl_v, 1)}$,\dots, $\phi_{3,1}^{(\bl_v, M^d)}$, 
$\phi_{4,1}^{(1)}, \dots, \phi_{4,1}^{(\lceil c_{13} \cdot M^{d^*}\rceil)}$ and 
$\bm{\phi}_{1,2}, \bm{\phi}_{2,2}, \phi_{3,2}^{(\bl_v)}$ 
$(v \in \{1, \dots, \binom{d+q}{d}\})$ without an error. Thus, it 
follows that 
\begin{align*}
\left|\bm{\hat{\phi}}_{1,2} - \bm{\hat{\phi}}_{2,2}\right| = \left|\bold{x}-\bm{\phi}_{2,2}\right| \leq 2a
\end{align*}
and 
\begin{align*}
\left|\hat{\phi}_{3, 2}^{(\bl_v)}\right| = \left|\phi_{3, 2}^{(\bl_v)}\right| \leq \|f\|_{C^{q}([-a,a]^d)}.
\end{align*}
Therefore, the input of $\hat{f}_p$ in \eqref{fp} is contained 
in the interval where \eqref{fpeq} holds. By choosing
\begin{align*}
B_{M,p} = \lceil \log_4\left(M^{2p}\right)\rceil
\end{align*}
we get
\begin{align*}
  &
  \left|\hat{f}_{\mathcal{P}_2}(\bold{x}) - 
T_{f,q,(C_{\mathcal{P}_2}(\bold{x}))_{\mathrm{left}}}(\bold{x})\right|
  =
  \left|\hat{\phi}_{1,3} - \phi_{1,3}\right|
   \\
& \leq c_{16} \cdot \left(\max\left\{2a, \|f\|_{C^{q}(\R^d)}\right\}\right)^{4(q+1)} \cdot \frac{1}{M^{2p}}, 
\end{align*}
where we have used $\bar{r}(p) \leq 1$. 
This together with Lemma~\ref{le1} and (\ref{se3eq*})
implies the first assertion of the lemma. 
\\
\\
In the \textit{last step of the proof} we bound |$\hat{f}_{\mathcal{P}_2}(\bold{x})|$ in case that $\bx \in \Rd$.
If $C_\bi \cap \M \neq \emptyset$ we know
\[
(\tilde{C}_{j,\bi})_{\mathrm{left}}
\in
   [-2a,2a]^d
   \quad
   (j \in \{1, \dots, M^d\}), 
\]
from which we can conclude
\begin{align*}
  \left|\hat{\phi}_{3, 1}^{(\bl, j)}\right| \leq
\|f\|_{C^q (\Rd)}
 \quad (j \in \{1, \dots,
  \lceil c_{13} \cdot M^{d^*} \rceil \})
\end{align*}
and 
\begin{align*}
  \left|\mathbf{\hat{\phi}}_{2,1}^{(j,s)}\right| \leq 2 \cdot a
  \quad (j \in \{1, \dots,
  \lceil c_{13} \cdot M^{d^*} \rceil \}, s \in \{1, \dots, d\}).
\end{align*}
Here we have used, that the value of $\hat{f}_{\textrm{ind}, C_{\bk}}$ lies (due to its construction in Lemma 5a)) in the interval $[0,1]$ and that
for fixed $\bx \in \Rd$ at most one of the values
$\hat{f}_{ind, C_\bi}(\bx)$ $(\bi \in \Z^d)$ is not equal to zero.
To bound the values of $\hat{\phi}_{3, 2}^{(\bl)}$ and $\hat{\phi}_{2,2}^{(j,s)}$ we consider the sums in \eqref{neur31} and \eqref{neur32}.  Due to the fact that all cubes 
$[\bm{\hat{\phi}}_{2,1}^{(j)},  \bm{\hat{\phi}}_{2,1}^{(j)} + \hat{\phi}_{4,1}^{(j)}\cdot \mathbf{1})$ are distinct for different $j \in \{1, \dots, \lceil c_{13} \cdot M^{d^*}\rceil\}$, those sums produce for at most one summand a value not equal to zero.  By construction of $\hat{f}_{\textrm{test}}$ in Lemma 5 this value, in turn,  is bounded
in absolute value
  by $|\hat{\phi}_{2,1}^{(j, s)}|$ or $|\hat{\phi}_{3,1}^{(\bl, j)}|$, respectively. 
This leads to 
\begin{align*}
\left|\hat{\phi}_{3, 2}^{(\bl)}\right| \leq \|f\|_{C^{q}(\R^d)}
\end{align*}
and 
\begin{align*}
\left|\hat{\phi}_{2,2}^{(s)}\right| \leq 2 \cdot a, \quad (s \in \{1, \dots, d\}).
\end{align*}
We conclude
\begin{eqnarray*}
\left|\hat{f}_{\mathcal{P}_2}(\bold{x})\right| &\leq &\left|\hat{f}_p\left(\bold{z},y_1, \dots, y_{\binom{d+q}{d}}\right) - p\left(\bold{z},y_1, \dots, y_{\binom{d+q}{d}}\right)\right|\\
&&+ \left|p\left(\bold{z},y_1, \dots, y_{\binom{d+q}{d}}\right)\right|\\
&\leq & 1 + \sum_{0 \leq \|\bl\|_1 \leq q} \frac{1}{\bl!} \cdot \|f\|_{C^{q}(\R^d)} \cdot \left(4a\right)^{\|\bl\|_1}\\
&
\leq&
1 +
\|f\|_{C^{q}(\R^d)}
\cdot
\left(
\sum_{l=0}^\infty \frac{(4a)^l}{l!}
\right)^d
\\
&
=&  1+ e^{4ad} \cdot \|f\|_{C^{q}(\R^d)}.
\end{eqnarray*}
\hfill $\Box$

\subsection{Key step 3 of the proof of Theorem~\ref{th2}: 
Approximation of $w_{\P_2}(\bold{x}) \cdot f(\bold{x})$ by neural 
networks }
In our key step 3 we construct a network that approximates
\[
w_{\P_2}(\bold{x}) \cdot f(\bold{x}),
\]
where
\begin{equation}
  \label{w_vb}
w_{\P_2}(\bold{x}) = \prod_{j=1}^d \left(1- 2 \cdot M^2 \cdot 
\left|(C_{\mathcal{P}_{2}}(\bold{x}))_{\mathrm{left}}^{(j)} + 
\frac{1}{2 \cdot M^2} - x^{(j)}\right|\right)_+
\end{equation}
is a linear tensorproduct B-spline
which takes its maximum value at the center of $C_{\P_{2}}(\bold{x})$, which
is nonzero in the inner part of $C_{\P_{2}}(\bold{x})$ and which
vanishes
outside of $C_{\P_{2}}(\bold{x})$. 

\begin{lemma}
\label{le6}
Let $\sigma: \R \to \R$ be the ReLU activation function $\sigma(x) = \max\{x,0\}$.  Let
                 $\M$ be a $d^*$-dimensional Lipschitz-manifold and let
                 $1 \leq a < \infty$ such that
$\M \subseteq [-a,a]^d$.  Let $p=q+s$ for some $q \in \N_0$, $s \in (0,1]$ and let $C>0$.
    Let $f: \Rd \to \R$ be a $(p,C)$-smooth function and let $w_{\P_2}$ be defined as in
\eqref{w_vb}. Let $M \in \N_0$ be such that 
    \begin{align*}
      M^{2p} \geq
      c_{17} \cdot \left(\max\left\{3a, \|f\|_{C^q(\R^d)}\right\}\right)^{4(q+1)}
    \end{align*}
    and
    \[
M^{2p} \geq c_{18} \cdot (2 \cdot a \cdot d)^{2p} \cdot C
    \]
hold.
Then there exists a network
\begin{align*}
\hat{f} \in \mathcal{F}\left(L, r\right)
\end{align*}
with
\begin{align*}
L=5+\lceil \log_4(M^{2p})\rceil \cdot \left(\lceil \log_2(\max\{q, d\}+1\})\rceil+1\right)
\end{align*}
and
\begin{align*}
  r=& 64 \cdot \binom{d+q}{d} \cdot d^2 \cdot (q+1) \cdot \lceil c_{13} \cdot
  M^{d^*} \rceil
\end{align*}
such that
\begin{align*}
&\left|\hat{f}(\bold{x}) - w_{\P_2}(\bold{x}) \cdot f(\bold{x})\right| \leq c_{19} \cdot \left(\max\left\{3a,  \|f\|_{C^q(\R^d)}\right\}\right)^{4(q+1)} \cdot \frac{1}{M^{2p}}
\end{align*}
holds for $\bold{x} \in \M$. 
\end{lemma}

In the proof of Lemma~\ref{le6} we adapt the arguments
in the proof of Lemma~7 in Supplement~A of Langer and Kohler (2020)
to the case that our input is contained in a $d^*$-dimensional
Lipschitz-manifold. To do this, we need the following two auxiliary results.

\begin{lemma}
\label{lea8}
Let $\sigma: \R \to \R$ be the ReLU activation function $\sigma(x) = \max\{x,0\}$.
Let
$\M$ be a $d^*$-dimensional Lipschitz-manifold and let
$1 \leq a < \infty$ such that
$\M \subseteq [-a,a]^d$. Let
$M \geq 4^{4d+1} \cdot d$. Let $\mathcal{P}_{2}$
be the partition defined in \eqref{partition} and let 
$w_{\P_2}(\bold{x})$ be 
the corresponding weight defined by \eqref{w_vb}. Then there exists a 
neural network
\begin{align*}
\hat{f}_{w_{\P_2}} \in \mathcal{F}\left(5+\lceil \log_4(M^{2p})\rceil 
\cdot \lceil \log_2(d)\rceil, r \right)
\end{align*}
with\begin{align*}
r=\max\left\{18d, 2d+d \cdot \lceil c_{13} \cdot M^{d^*}\rceil  \cdot 2 \cdot (2+2d)\right\}
\end{align*}
such that
\begin{align*}
\left|\hat{f}_{w_{\P_2}}(\bold{x}) - w_{\P_2}(\bold{x})\right| \leq 4^{4d+1} \cdot d \cdot \frac{1}{M^{2p}}
\end{align*}
for $\bold{x} \in
\left(
\bigcup_{\bk \in \Z^d} (D_{\bk})_{1/M^{2p+2}}^0
\right)
\cap \M$
and 
\begin{align*}
|\hat{f}_{w_{\P_2}}(\bold{x})| \leq 2
\end{align*}
for $\bold{x} \in \M$.
\end{lemma}

\noindent
    {\bf Proof.}
    The proof follows by a slight modification from the proof of Lemma 9 in the Supplement A of Kohler and Langer (2020).  A complete proof is given in the appendix.
    
  \hfill $\Box$

\begin{lemma}
\label{lea9}
Let $\sigma: \R \to \R$ be the ReLU activation function $\sigma(x) = \max\{x,0\}$.
Let $\M$ be a $d^*$-dimensional Lipschitz-manifold and let
$1 \leq a < \infty$ such that
$\M \subseteq [-a,a]^d$. Let
$\mathcal{P}_{1}$ and $\mathcal{P}_{2}$
be the partitions defined in \eqref{partition} and let $M \in \N$. Then there exists a neural network 
\begin{align*}
\hat{f}_{\mathrm{check}, \mathcal{P}_{2}} \in \mathcal{F}\left(5, 2d + 
(4d^2+4d) \cdot \lceil c_{13} \cdot M^{d^*}\rceil \right)
\end{align*}
satisfying
\begin{align*}
  \hat{f}_{\mathrm{check}, \mathcal{P}_{2}}(\bold{x}) = \mathds{1}_{
    \bigcup_{\bi \in \Z^d}
    D_{\bk} \setminus (D_{\bk})_{1/M^{2p+2}}^0
}(\bold{x})
\end{align*}
for $\bold{x} \in \M \setminus
\left( \bigcup_{\bk \in \Z^d} (D_{\bk})_{1/M^{2p+2}}^0 \textbackslash (D_{\bk})_{2/M^{2p+2}}^0
\right)$ and 
\begin{align*}
\hat{f}_{\mathrm{check}, \mathcal{P}_{2}}(\bold{x}) \in [0,1]
\end{align*}
for $\bold{x} \in \M$. 
\end{lemma}

\noindent
    {\bf Proof.}
Throughout the proof we assume that $\bi \in \Z^d$ satisfies $C_{\P_1}(\bold{x}) = C_{\bi}$. 
In oder to compute $\hat{f}_{\mathrm{check}, \mathcal{P}_2}$ we use a 
two-scale approximation
defined as follows: In the first part of the network we check whether $\bold{x}
\in \M$ is contained in
\begin{align*}
\bigcup_{\bk \in \Z^d} C_{\bk}\setminus (C_{\bk})_{1/M^{2p+2}}^0.
\end{align*}
Therefore, our network approximates in the first two hidden layers
for $\bx \in \M$
the function
\begin{align*}
  f_1(\bold{x}) &= \mathds{1}_{\bigcup_{
\bk \in \Z^d : C_\bk \cap \M \neq \emptyset
    } C_{\bk}\setminus (C_{\bk})_{1/M^{2p+2}}^0}(\bold{x})=1-\sum_{
\bk \in \Z^d : C_\bk \cap \M \neq \emptyset
} \mathds{1}_{(C_{\bk})_{1/M^{2p+2}}^0}(\bold{x})
\end{align*}
by 
\begin{align*}
  \hat{f}_1(\bold{x})= 1-\sum_{
\bk \in \Z^d : C_\bk \cap \M \neq \emptyset
  } \hat{f}_{\mathrm{ind}, (C_{\bk})_{1/M^{2p+2}}^0}(\bold{x}),
\end{align*}
where $\hat{f}_{\mathrm{ind}, (C_{\bk})_{1/M^{2p+2}}^0}$ are the 
networks of
Lemma~\ref{le5}a), which need $2d$ neurons per layer, respectively. To 
approximate the indicator functions on the partition $\P_2$ only for 
the cubes
$D_{\bk} \subset C_{\P_1}(\bold{x}) = C_{\bi}$, we further need to
compute the positions of $(\tilde{C}_{j, \bi})_{\mathrm{left}}$ 
$(j \in \{1, \dots, \lceil c_{13} \cdot M^{d^*}\rceil\})$.  This can 
be done as 
described by the networks $\bm{\hat{\phi}}_{2,1}^{(j)}$ in the proof of 
Lemma~\ref{le2} with
$d \cdot \lceil c_{13} \cdot M^{d^*} \rceil \cdot 2d$ neurons.  The length of 
each cube $(C_{j, \bi})_{\textrm{left}}$ $(j \in \{1, \dots, \lceil c_{13} \cdot M^{d^*}\rceil\}$ is computed by $\hat{\phi}_{4,1}^{(j)}$ as in the proof of Lemma \ref{le2} with $\lceil c_{13} \cdot M^{d^*} \rceil \cdot 2d$ neurons per layer. To 
shift the value of $\bold{x}$ 
in the next hidden layers we further apply the network 
$\hat{f}^2_{\mathrm{id}}$,  which needs $2d$ neurons per layer. 
Analogous to (\ref{Aj}) we can describe the cubes 
$(\tilde{C}_{j, \bi})_{1/M^{2p+2}}^0$ $(j \in \{1, \dots, N_\bi\})$
by 
\begin{align*}
&(\mathcal{A}^{(j)})_{1/M^{2p+2}}^0 = \left\{\bold{x} \in \Rd: -x^{(k)} + \phi_{2,1}^{(j,k)} +\frac{1}{M^{2p+2}}\leq 0 \right.\\
 & \left.  \hspace{3cm} \text{and} \ 
 x^{(k)} - \phi_{2,1}^{(j,k)} - 
\phi_{4,1}^{(j)} +\frac{1}{M^{2p+2}} < 0 \ \text{for all} \ k \in \{1, 
\dots, d\}\right\}.
\end{align*}
Then, for $\bx \in \M$, the function 
\begin{align*}
  f_2(\bold{x}) = \mathds{1}_{\bigcup_{
\bk \in Z^d: C_\bk \cap \M \neq \emptyset, j \leq N_\bk
    } \tilde{C}_{j,\bk} \setminus (\tilde{C}_{j,\bk})_{1/M^{2p+2}}^0}(\bold{x})=
  1-\sum_{
\bk \in Z^d: C_\bk \cap \M \neq \emptyset, j \leq N_\bk
  } \mathds{1}_{(\tilde{C}_{j,\bk})_{1/M^{2p+2}}^0}(\bold{x})
\end{align*}
can be approximated by 
\begin{align*}
\hat{f}_2(\bold{x}) &= 
1 - \sum_{j \in \{1, \dots, \lceil c_{13} \cdot M^{d^*}\rceil 
\}}\hat{f}_{\mathrm{test}}\left(\hat{f}_{\mathrm{id}}^2(\bold{x}), 
\bm{\hat{\phi}}_{2,1}^{(j)} +\frac{1}{M^{2p+2}}\cdot \mathbf{1}, \right.\\
& \hspace*{5cm} \left.\bm{\hat{\phi}}_{2,1}^{(j)} +
\hat{\phi}_{4,1}^{(j)} \cdot \mathbf{1}-\frac{1}{M^{2p+2}}\cdot \mathbf{1}, 
1\right), 
\end{align*}
\noindent where $\hat{f}_{\mathrm{test}}$ is the network of Lemma~\ref{le5}b), 
which needs $2$ hidden layers and $2 \cdot (2d+2)$ neurons per layer.
Here, for any $\bx \in \Rd$ at most one of the terms in the
sum in the definition of $\hat{f}_2(\bold{x})$ is not equal to zero,
and $\hat{\phi}_{4,1}^{(j)}$ is equal to zero for $j > N_\bi$.
\\
Combining the networks $\hat{f}_1$ and $\hat{f}_2$ and using the 
characteristics of ReLU activation function that is zero in case of 
negative input, finally let us approximate 
\[
\mathds{1}_{
    \bigcup_{\bk \in \Z^d}
    D_{\bk} \setminus (D_{\bk})_{1/M^{2p+2}}^0
}(\bold{x})
\]
for
$\bx \in \M$
by
\begin{align*}
\hat{f}_{\mathrm{check}, \mathcal{P}_{2}}(\bold{x}) &= 
1-\sigma\left(1-\hat{f}_2(\bold{x}) - 
\hat{f}_{\mathrm{id}}^2\left(\hat{f}_1(\bold{x})\right)\right)\\
&= 1-\sigma\left(\sum_{j \in \{1, \dots, \lceil c_{13} \cdot 
M^{d^*}\rceil 
\}}\hat{f}_{\mathrm{test}}\left(\hat{f}_{\mathrm{id}}^2(\bold{x}), 
\bm{\hat{\phi}}_{2,1}^{(j)} +\frac{1}{M^{2p+2}}\cdot \mathbf{1}, 
\right.\right.\\
& \hspace*{2cm} \left. \bm{\hat{\phi}}_{2,1}^{(j)} +\hat{\phi}_{4,1}^{(j)} \cdot \mathbf{1}-\frac{1}{M^{2p+2}}\cdot \mathbf{1}, 1\right) \\
& \left. \hspace*{1.2cm} - \hat{f}_{\mathrm{id}}^2\left(1-\sum_{\bk 
\in \Z^d: C_{\bk} \cap \M \neq \emptyset}\hat{f}_{\mathrm{ind}, 
(C_{\bk})_{1/M^{2p+2}}^0}(\bold{x})\right)\right). 
\end{align*}
Now it is easy to see that our whole network is contained in the 
network class
\begin{align*}
\mathcal{F}(5, r)
\end{align*}
with
\begin{align*}
  r&=\max\{2d+d \cdot \lceil c_{13} \cdot M^{d^*} \rceil
  \cdot 2d+
  \lceil c_{13} \cdot M^{d^*} \rceil\cdot 2d, \lceil c_{13} \cdot M^{d^*} \rceil \cdot 2 \cdot (2+2d)+2\}\\
&\leq 2d + (4d^2+4d) \cdot \lceil c_{13} \cdot M^{d^*} \rceil. 
\end{align*}
As in the proof of Lemma 10 in Kohler and Langer (2020) it can be shown that 
\begin{align}
\label{fcheckeq1}
\hat{f}_{\mathrm{check}, \P_2}(\bold{x})= \mathds{1}_{\bigcup_{\bk \in 
\Z^d} D_{\bk} \setminus (D_{\bk})_{1/M^{2p+2}}^0}(\bold{x})
\end{align}
holds for $\bold{x} \in \M \setminus
\left( \bigcup_{\bk \in \Z^d} (D_{\bk})_{1/M^{2p+2}}^0 \textbackslash (D_{\bk})_{2/M^{2p+2}}^0
\right)$.  This part of the proof is given in the appendix.
    \hfill $\Box$
\\  
\\
  \noindent
      {\bf Proof of Lemma~\ref{le6}.} Using the networks $\hat{f}_{\mathcal{P}_{2}}$ of Lemma \ref{le2},  $\hat{f}_{\mathrm{check}, \mathcal{P}_{2}}$ of Lemma \ref{lea9} and $\hat{f}_{w_{\P_2}}$ of Lemma \ref{lea8}, this proof follows directly from the proof of Lemma 7 in the Supplement A of Kohler and Langer (2020). A complete proof is given in the appendix.

\subsection{Key step~4 of the proof of Theorem~\ref{th2}: 
Applying $\hat{f}$ to slightly shifted partitions}
Finally, we will use a finite sum of the networks of Lemma~\ref{le6}, where $\P_2$ is substituted by a slightly shifted version of $\P_2$, respectively, to approximates $f(\bx)$ in supremum norm and to show Theorem \ref{th2}.  
\\
\\
%
\noindent
{\bf Proof of Theorem~\ref{th2}.}
The proof follows as a slight modification from the proof of Theorem 2 in Kohler and Langer (2020), where we use the networks $\hat{f}$ of Lemma \ref{le6}. A complete proof is given in the appendix.
\hfill $\Box$

{\footnotesize{
}}
\newpage
\section*{Appendix}
{\bf An auxiliary result for the proof of Lemma \ref{lea8}.}
In the proof of Lemma \ref{lea8} we will need the following
auxiliary result.

\begin{lemma}
\label{nle1}
Let $\sigma: \R \to \R$ be the ReLU activation function $\sigma(x) = 
\max\{x,0\}$. 
Then for any $R \in \N$ and any $b \geq 1$ a neural network 
\begin{align*}
\hat{f}_{\mathrm{mult}, d} \in \mathcal{F}(R \cdot \lceil \log_2(d) 
\rceil, 18d)
\end{align*}
exists such that
\begin{align*}
\left|\hat{f}_{\mathrm{mult}, d}(\bold{x}) - 
\prod_{i=1}^dx^{(i)}\right| \leq 4^{4d+1}\cdot b^{4d} \cdot d \cdot 
4^{-R}
\end{align*}
holds for all $\bold{x} \in [-b,b]^d$. 
\end{lemma}

\noindent
{\bf Proof.}
See Lemma~8 in Supplement~A of Kohler and Langer (2020).
\hfill $\Box$\\

\noindent
{\bf Proof of Lemma \ref{lea8}.}
 The first four hidden layers of $\hat{f}_{w_{\P_2}}$ compute
for $\bx \in \M$
    the value of 
\begin{align*}
(C_{\mathcal{P}_{2}}(\bold{x}))_{\mathrm{left}}
\end{align*}
and shift the value of $\bold{x}$ in the next hidden layer, 
respectively. 
This can be done as described in $\bm{\hat{\phi}}_{1,2}$ and 
$\bm{\hat{\phi}}_{2,2}$ in the proof of Lemma~\ref{le2} with $2d+d 
\cdot
\lceil c_{13} \cdot M^{d^*} \rceil
\cdot 2 \cdot (2+2d)$ neurons per layer.
The fifth hidden layer then computes the functions
\begin{eqnarray*}
&&\left(1-2 \cdot M^2 \cdot 
\left|(C_{\mathcal{P}_{2}}(\bold{x}))_{\mathrm{left}}^{(j)} + 
\frac{1}{2 \cdot M^2} -x^{(j)}\right| \right)_+\\
  &&=
  \left(
  2 \cdot M^2 \cdot
  \left(
x^{(j)} - (C_{\mathcal{P}_{2}}(\bold{x}))_{\mathrm{left}}^{(j)} 
  \right)
  \right)_+
  \\
  &&
  \quad
  -
  2 \cdot \left(
  2 \cdot M^2 \cdot
  \left(
  x^{(j)} - (C_{\mathcal{P}_{2}}(\bold{x}))_{\mathrm{left}}^{(j)}
  - \frac{1}{2 \cdot M^2}
  \right)
  \right)_+
 \\
  &&
  \quad
  +
  \left(
  2 \cdot M^2 \cdot
  \left(
  x^{(j)} - (C_{\mathcal{P}_{2}}(\bold{x}))_{\mathrm{left}}^{(j)}
  - \frac{1}{M^2}
  \right)
  \right)_+, \quad j \in \{1, \dots, d\},
\end{eqnarray*}
 using the networks
\begin{align*}
  \hat{f}_{w_{{\P_2},j}}(\bold{x}) &= \sigma\left(
  2 \cdot M^2 \cdot
  \left(
\hat{\phi}_{1,2}^{(j)} - \bm{\hat{\phi}}_{2,2}^{(j)} 
  \right)
  \right)\\
  & \quad -2 \cdot \sigma\left(
2 \cdot M^2 \cdot
  \left(
  \hat{\phi}_{1,2}^{(j)} - \hat{\phi}_{2,2}^{(j)}
  - \frac{1}{2 \cdot M^2}
  \right) 
  \right)\\
  & \quad + \sigma\left(
  2 \cdot M^2 \cdot
  \left(
  \hat{\phi}_{1,2}^{(j)} - \hat{\phi}_{2,2}^{(j)}
  - \frac{1}{M^2}
  \right)
  \right),  \quad j \in \{1, \dots, d\},
\end{align*}
with $3d$ neurons. The product \eqref{w_vb} of $w_{\P_2,j}(\bold{x})$ $(j \in \{1, \dots, d\})$
can then be computed by the network $\hat{f}_{\mathrm{mult},d}$ 
of Lemma~\ref{nle1}, where we choose $x^{(j)} = 
\hat{f}_{w_{{\P_2},j}}(\bold{x})$. 
Finally we set
\begin{align*}
\hat{f}_{w_{\P_2}}(\bold{x}) = \hat{f}_{\mathrm{mult}, 
d}\left(\hat{f}_{w_{{\P_2},1}}(\bold{x}), \dots, 
\hat{f}_{w_{{\P_2},d}}(\bold{x})\right).
\end{align*}
By choosing $R= \lceil \log_4(M^{2p})\rceil$ in Lemma~\ref{nle1}, this 
network lies in the class
\begin{align*}
  \mathcal{F}\left(4+1+\lceil \log_4(M^{2p})\rceil \cdot \lceil\log_2(d)\rceil, \max\left\{18d,2d+ d \cdot \lceil c_{13} \cdot M^{d^*} \rceil
  \cdot 2 \cdot (2+2d), 3d\right\}\right),
\end{align*}
and
according to Lemma~\ref{nle1} (where we set $b=1$) it
approximates $w_{\P_2}(\bold{x})$ with an error of size
\begin{align*}
4^{4d+1} \cdot d \cdot \frac{1}{M^{2p}}
\end{align*}
in case that $\bold{x} \in \M$ is contained in
$\bigcup_{\bk \in \Z^d} (D_{\bk})_{1/M^{2p+2}}^0$. Since 
$|\hat{f}_{w_{\P_2},j}(\bold{x})| \leq 1$ for \linebreak $j \in \{1, 
\dots, d\}$ 
we can bound the value of the network using the triangle inequality by
\begin{align*}
|\hat{f}_{w_{\P_2}}(\bold{x})| \leq \left|\hat{f}_{w_{\P_2}}(\bold{x}) - \prod_{j=1}^d \hat{f}_{w_{{\P_2},j}}(\bold{x})\right| + \left|\prod_{j=1}^d \hat{f}_{w_{{\P_2},j}}(\bold{x})\right| \leq 2
\end{align*}
for $\bold{x} \in \M$, where we have used that
\begin{align*}
M^{2p} \geq 4^{4d+1} \cdot d.
\end{align*}

    \hfill $\Box$
    \\
\noindent
{\bf Network accuracy of $\hat{f}_{\mathrm{check}, \mathcal{P}_{2}}$ in Lemma \ref{lea9}.}
\\
{\bf Proof of \eqref{fcheckeq1}.}
We distinguish between \textit{three} cases. In our first case we assume that 
\begin{align*}
\bx \in \M \ \text{and} \ \bold{x} \notin \bigcup_{\bk \in \Z^d} 
(C_{\bk})_{1/M^{2p+2}}^0, 
\end{align*}
which  implies that
\begin{align*}
\bold{x} \notin \bigcup_{\bk \in \Z^d} (D_{\bk})_{1/M^{2p+2}}^0.
\end{align*}
In this case we get from 
Lemma~\ref{le5}  that
$\hat{f}_1(\bold{x})=1$ from which we can conclude
\begin{align*}
&1-\hat{f}_2(\bold{x}) - 
\hat{f}_{\mathrm{id}}^2\left(\hat{f}_1(\bold{x})\right) \\
&=\sum_{j \in \{1, \dots, \lceil c_{13} \cdot M^{d^*} \rceil 
\}}\hat{f}_{\mathrm{test}}\left(\hat{f}_{\mathrm{id}}^2(\bold{x}), 
\bm{\hat{\phi}}^{(j)}_{2,1} +\frac{1}{M^{2p+2}}\cdot \mathbf{1}, 
\bm{\hat{\phi}}^{(j)}_{2,1} +\hat{\phi}_{4,1}^{(j)} \cdot 
\mathbf{1}-\frac{1}{M^{2p+2}}\cdot \mathbf{1}, 1\right) -1\\
&\leq 0.
\end{align*}
Here we have used that each $\hat{f}_{\mathrm{test}}$ is contained in 
$[0,1]$ (according to its construction in Lemma 5b)) and that at most one $\hat{f}_{\mathrm{test}}$ in the sum is 
larger than $0$. Finally we get
\[
\hat{f}_{\mathrm{check}, \P_2}(\bold{x})=1-0=1=\mathds{1}_{
    \bigcup_{\bk \in \Z^d}
    D_{\bk} \setminus (D_{\bk})_{1/M^{2p+2}}^0
}(\bold{x}).
\]
In our second case we assume that 
\begin{align}
\label{eq300}
\bold{x} \in \M \cap \left(\bigcup_{\bk \in \Z^d} (C_{\bk})_{1/M^{2p+2}}^0\right).
\end{align}
and 
\begin{align*}
\bold{x} \in \bigcup_{\bk \in \Z^d} (D_{\bk})_{2/M^{2p+2}}^0.
\end{align*} 

 Then we have $\hat{\phi}_{2,1}^{(j)} = (\tilde{C}_{j, 
   \bi})_{\mathrm{left}}$ and
 $\hat{\phi}_{4,1}^{(j)} = (1/M^2) \cdot \mathds{1}_{\{j \leq N_\bi \}}$.
  Furthermore, we can conclude that
 \begin{align*}
   (\mathcal{A}^{(j)})_{1/M^{2p+2}}^0 &= \Bigg\{\bold{x} \in \Rd: -\hat{\phi}_{1,1}^{(k)} + \hat{\phi}_{2,1}^{(j,k)} +\frac{1}{M^{2p+2}}\leq 0 \\
 & \quad \quad \quad \text{and} \ 
 \hat{\phi}_{1,1}^{(k)} - 
\hat{\phi}_{2,1}^{(j, k)}  - \hat{\phi}_{4,1}^{(j)} +\frac{1}{M^{2p+2}} < 0 \\
&  \quad \quad \quad \text{for all} \ k \in \{1, \dots, d\}\Bigg\}\\
 &= (\tilde{C}_{j, \bi})_{1/M^{2p+2}}^0 
\end{align*}
 for $j \in \{1, \dots, N_\bi\}$. Since we only have to show our assumption for  
\begin{align*}
\bold{x} \notin \bigcup_{\bk \in \Z^d} (D_{\bk})_{1/M^{2p+2}}^0 \textbackslash (D_{\bk})_{2/M^{2p+2}}^0, 
\end{align*} 
%
%
%
%
we can conclude by Lemma~\ref{le5}  that
\begin{eqnarray*}
  &&
  \hat{f}_{\mathrm{test}}\Big(\bm{\hat{\phi}}_{1,1}, \bm{\hat{\phi}}^{(j)}_{2,1} 
+\frac{1}{M^{2p+2}} \cdot \mathbf{1}, 
\bm{\hat{\phi}}_{2,1}^{(j)} +\hat{\phi}_{4,1}^{(j)} \cdot \mathbf{1} -\frac{1}{M^{2p+2}}\cdot \mathbf{1}, 1 \Big)
=
\mathds{1}_{ (\tilde{C}_{j,\bi})_{1/M^{2p+2}}^0}(\bold{x}) \cdot \mathds{1}_{\{j \leq N_\bi \}}
\end{eqnarray*}
for all $j \in \{1, \dots, \lceil c_{13} \cdot M^{d^*}\rceil \}$.
This implies
\[
\hat{f}_2(\bold{x}) = f_2(\bold{x}).
\]
Since
\begin{align*}
\bold{x} \in \bigcup_{\bk \in \Z^d} (D_{\bk})_{2/M^{2p+2}}^0
\end{align*} 
we can further conclude that
\begin{align*}
\bold{x} \in \bigcup_{\bk \in \Z^d} (C_{\bk})_{2/M^{2p+2}}^0
\end{align*}
and it follows by Lemma~\ref{le5} that
\begin{align*}
\hat{f}_1(\bold{x}) = f_1(\bold{x})=0.
\end{align*}
Thus, we have
\begin{align*}
1-\hat{f}_2(\bold{x}) -\hat{f}_{\mathrm{id}}^2(\hat{f}_1(\bold{x})) = 
1-f_2(\bold{x}) = 1-0 = 1
\end{align*}
and 
\begin{align*}
\hat{f}_{\mathrm{check}, \P_2}(\bold{x}) = 1-1=0= 
\mathds{1}_{\bigcup_{\bk \in \Z^d} D_{\bk} \setminus 
(D_{\bk})_{1/M^{2p+2}}^0}(\bold{x}).
\end{align*}
In our third case we assume \eqref{eq300}, but 
\begin{align*}
\bold{x} \in \bigcup_{\bk \in \Z^d} (D_{\bk}) \textbackslash (D_{\bk})^0_{1/M^{2p+2}}, 
\end{align*}
which means that
\begin{align*}
\bold{x} \notin \bigcup_{\bk \in \Z^d} (D_{\bk})^0_{1/M^{2p+2}}.
\end{align*}
In this case the approximation $\hat{f}_1(\bold{x})$ is not exact. By the definition of the networks in Lemma~\ref{le5}a),  we can conclude that all values of $\hat{f}_{\mathrm{ind}, 
(C_{\bk})_{1/M^{2p+2}}^0}$ in the definition of $\hat{f}_1$ are 
contained in $[0,1]$
and that at most one of them is greater than zero.
Thus, we have
 \begin{align*}
 \hat{f}_1(\bold{x}) \in [0,1].
 \end{align*}
 Since \eqref{eq300} holds we further have
 \begin{align*}
 \hat{f}_2(\bold{x}) = f_2(\bold{x})
 \end{align*}
 as shown in the second case. Summarizing, we can conclude that
\begin{align*}
1-\hat{f}_2(\bold{x}) - \hat{f}_{\mathrm{id}}^2(\hat{f}_1(\bold{x})) &= 
\sum_{j \in \{1, \dots, \lceil c_{13} \cdot M^{d^*}\rceil \}} 
\mathds{1}_{(\tilde{C}_{j,\bi})_{1/M^{2p+2}}^0}(\bold{x}) - 
\hat{f}_{\mathrm{id}}^2(\hat{f}_1(\bold{x}))\\
&\leq 0-0 = 0.
\end{align*}
This implies
\begin{align*}
\hat{f}_{\mathrm{check}, \P_2}(\bold{x}) = 1-0 =1 
=\mathds{1}_{\bigcup_{\bk \in \Z^d} D_{\bk} \setminus 
(D_{\bk})_{1/M^{2p+2}}^0}(\bold{x}).
\end{align*}
By construction of the network we have
$\hat{f}_1(\bx) \geq 0$ and $\hat{f}_2(\bx) \geq 0$
$(\bx \in \Rd)$, hence
\begin{align*}
\hat{f}_{\mathrm{check}, \mathcal{P}_{2}}(\bold{x}) \in [0,1]
\end{align*}
holds for $\bold{x} \in \M$. 
\hfill $\Box$
\\
\noindent
{\bf An auxiliary result for the proof of Lemma \ref{le6}.}
In the proof of Lemma \ref{le6} we will need the following auxiliary result.

\begin{lemma}\label{le3}
Let $\sigma: \mathbb{R} \to \R$ be the ReLU activation function 
$\sigma(x) = \max\{x,0\}$. 
Then for any $R \in \N$ and any $b \geq 1$ a neural network
\begin{equation*}
\hat{f}_{\mathrm{mult}} \in \mathcal{F}(R,18)
\end{equation*}
 exists such that
\begin{equation*}
|\hat{f}_{\mathrm{mult}}(x,y) - x \cdot y| \leq 2 \cdot b^2 \cdot 
4^{-R}
\end{equation*}
holds for all $x, y \in [-b, b]$.
\end{lemma}
\noindent
{\bf Proof.}
See Lemma~4 in Supplement~A of Kohler and Langer (2020).
\hfill $\Box$
\\    

\noindent
{\bf Proof of Lemma \ref{le6}.}
        Let $\hat{f}_{\mathcal{P}_{2}}$ be the network of 
Lemma~\ref{le2}
  and let $\hat{f}_{\mathrm{check}, \mathcal{P}_{2}}$ be the network 
of Lemma~\ref{lea9}. 
  By successively applying $\hat{f}_{\mathrm{id}}$ to the output of 
  one of these networks, we can achieve that both networks have the 
  same number of hidden layers, i.e.,
  \begin{align*}
  L=4+\max\left\{\lceil \log_4(M^{2p})\rceil \cdot \lceil\log_2(\max\{q+1, 2\})\rceil, 1\right\}.
  \end{align*}
  We set
\begin{eqnarray*}
  \hat{f}_{\mathcal{P}_{2}, \mathrm{true}}(\bold{x}) &=& 
\sigma\left(\hat{f}_{\mathcal{P}_{2}}(\bold{x}) - B_{\mathrm{true}} 
\cdot \hat{f}_{\mathrm{check}, \mathcal{P}_{2}}(\bold{x})\right)\\
  &&
  - \sigma\left(-\hat{f}_{\mathcal{P}_{2}}(\bold{x}) - 
B_{\mathrm{true}} \cdot \hat{f}_{\mathrm{check}, 
\mathcal{P}_{2}}(\bold{x})\right),
\end{eqnarray*}
where 
\begin{align*}
B_{\mathrm{true}} = 2 \cdot e^{4ad} \cdot \max 
\left\{\|f\|_{C^q(\R^d)},1\right\}.
\end{align*}
This network is contained in den network class $\mathcal{F}(L,r)$ with
\begin{align*}
L= 5+\lceil \log_4(M^{2p})\rceil \cdot \lceil \log_2(\max\{q+1, 2\})\rceil
\end{align*}
and
\begin{align*}
r=&
\max\left\{\left(\binom{d+q}{d} + d\right) \cdot \lceil c_{13} \cdot 
M^{d^*}\rceil \cdot 2 \cdot (2+2d)+2d, 18 \cdot (q+1) \cdot 
\binom{d+q}{d}\right\}\\
&+2d+(4d^2+4d) \cdot \lceil c_{13} \cdot M^{d^*}\rceil.
\end{align*}
Due to the fact that the value of $\hat{f}_{\mathcal{P}_{2}}$ is 
bounded by $B_{\mathrm{true}}$ according to Lemma~\ref{le2} and that 
$\hat{f}_{\mathrm{check}, \mathcal{P}_{2}}(\bold{x})$ is $1$ in case 
that $\bold{x} \in \M$ lies in
\begin{align}
\label{noset100}
\bigcup_{\bi \in \Z^d}
    D_{\bi} \setminus (D_{\bi})_{1/M^{2p+2}}^0,
\end{align} 
the properties of the ReLU activation function imply that the value 
of 
$\hat{f}_{\mathcal{P}_{2}, \mathrm{true}}(\bold{x})$ is zero in case 
that $\bold{x}$ is contained in \eqref{noset100}. 
Let $\hat{f}_{w_{\P_2}}$ be the network of Lemma~\ref{lea8}. 
To multiply the network $\hat{f}_{\mathcal{P}_{2},\mathrm{true}}$ by 
$\hat{f}_{w_{\P_2}}$ we use the network 
\begin{align*}
\hat{f}_{\mathrm{mult}} \in \mathcal{F}(\lceil\log_4(M^{2p})\rceil, 
18)
\end{align*}
of Lemma~\ref{le3}, which satisfies 
\begin{align}
\label{appfmult}
\left|\hat{f}_{\mathrm{mult}}(x,y) - xy\right| \leq 8 \cdot 
\left(\max\left\{\|f\|_{\infty} ,1\right\}\right)^2 \cdot 
\frac{1}{M^{2p}}
\end{align}
for all $x,y$ contained in
\begin{align*}
&\left[-2 \cdot \max\left\{\|f\|_{\infty} ,1\right\},  2\cdot \max\left\{\|f\|_{\infty},1\right\}\right].
\end{align*}
Here we have chosen $R=\lceil\log_4(M^{2p})\rceil$ in 
Lemma~\ref{le3}. 
\\
\\
By successively applying $\hat{f}_{\mathrm{id}}$ to the outputs of the 
networks $\hat{f}_{w_{\P_2}}$ and $\hat{f}_{\mathcal{P}_{2}, 
\mathrm{true}}$, we can synchronize their depths such that both 
networks have 
\begin{align*}
L=5+\lceil \log_4(M^{2p})\rceil \cdot \left(\lceil \log_2(\max\{q, d\}+1\})\rceil\right)
\end{align*}
hidden layers. 
\\
\\
The final network is given by
\begin{align*}
\hat{f}(\bold{x}) = 
\hat{f}_{\mathrm{mult}}\left(\hat{f}_{w_{\P_2}}(\bold{x}), 
\hat{f}_{\mathcal{P}_{2}, \mathrm{true}}(\bold{x})\right)
\end{align*}
and the network is contained in the network class $\mathcal{F}(L,r)$ with
\begin{align*}
L=5+\lceil \log_4(M^{2p})\rceil \cdot \left(\lceil \log_2(\max\{q, d\}+1\})\rceil+1\right) 
\end{align*}
and
\begin{align*}
r=&\max\left\{\left(\binom{d+q}{d} + d\right) \cdot \lceil c_{13} \cdot M^{d^*}\rceil \cdot 2 \cdot (2+2d)+2d, 18 \cdot (q+1) \cdot \binom{d+q}{d}\right\}\\
&+2d+(4d^2+4d) \cdot \lceil c_{13} \cdot M^{d^*}\rceil+\max\left\{18d, 2d+d \cdot \lceil c_{13} \cdot M^{d^*}\rceil \cdot 2 \cdot (2+2d)\right\}\\
\leq & 64 \cdot \binom{d+q}{d} \cdot d^2 \cdot (q+1) \cdot \lceil c_{13} \cdot M^{d^*}\rceil .
\end{align*}
In case that
\begin{align*}
\bold{x} \in \M \ \text{and} \ \bold{x} \in \bigcup_{\bk \in \Z^d} 
\left(D_{\bk}\right)_{2/M^{2p+2}}^0,
\end{align*}
the value of $\bold{x}$ is neither contained in
\begin{align}
\label{noset1}
\bigcup_{\bk \in \Z^d}
D_{\bk} \setminus (D_{\bk})^0_{1/M^{2p+2}}
\end{align}
nor contained in
\begin{align}
\label{noset2}
\bigcup_{\bk \in \Z^d}
(D_{\bk})^0_{1/M^{2p+2}} \setminus (D_{\bk})^0_{2/M^{2p+2}}
.
\end{align}
Thus, the network $\hat{f}_{w_{\P_2}}(\bold{x})$ approximates 
$w_{\P_2}(\bold{x})$ according to Lemma~\ref{lea8} with an error of 
size
\begin{align}
\label{appwv}
4^{d+1} \cdot d \cdot \frac{1}{M^{2p}}
\end{align}
and $\hat{f}_{\mathcal{P}_{2}}(\bold{x})$ approximates $f(\bold{x})$ according to
Lemma~\ref{le2} with an error of size
\begin{align}
\label{200}
c_{20} \cdot \left(\max\left\{2a, \|f\|_{C^q(\R^d)}\right\}\right)^{4(q+1)} \cdot \frac{1}{M^{2p}}.
\end{align}
Since $\hat{f}_{\mathrm{check}, \mathcal{P}_{2}}(\bold{x}) =0$, we 
have
\begin{align*}
  \hat{f}_{\mathcal{P}_{2}, \mathrm{true}}(\bold{x})
  =
  \sigma (\hat{f}_{\mathcal{P}_{2}}(\bold{x})) - \sigma( - \hat{f}_{\mathcal{P}_{2}}(\bold{x}))
  = \hat{f}_{\mathcal{P}_{2}}(\bold{x}).
\end{align*}
Since $M^{2p} \geq 4^{d+1} \cdot d$, we can bound the 
value of $\hat{f}_{w_{\P_2}}(\bold{x})$ using the triangle inequality 
by
\begin{align*}
|\hat{f}_{w_{\P_2}}(\bold{x})| \leq |\hat{f}_{w_{\P_2}}(\bold{x}) - w_{\P_2}(\bold{x})| + |w_{\P_2}(\bold{x})| \leq 2.
\end{align*}
Furthermore, we can bound 
\begin{align*}
  |\hat{f}_{\mathcal{P}_{2}}(\bold{x})| \leq |\hat{f}_{\mathcal{P}_{2}}(\bold{x})-f(\bold{x})| + |f(\bold{x})| \leq 2 \cdot
\max \{ \|f\|_{\infty},1 \}, 
\end{align*}
where we used $M^{2p} \geq c_{33} \cdot \left(\max\left\{2a, 
\|f\|_{C^q(\R^d)}\right\}\right)^{4(q+1)}$. Thus, the values of both networks are 
contained in the interval, where \eqref{appfmult} holds. 
Using the triangle inequality, this implies
\begin{align*}
&\left|\hat{f}_{\mathrm{mult}}\left(\hat{f}_{w_{\P_2}}(\bold{x}), 
\hat{f}_{\mathcal{P}_{2}, \mathrm{true}}(\bold{x})\right) - 
w_{\P_2}(\bold{x}) \cdot f(\bold{x})\right|\\
& \leq 
\left|\hat{f}_{\mathrm{mult}}\left(\hat{f}_{w_{\P_2}}(\bold{x}), 
\hat{f}_{\mathcal{P}_{2}, \mathrm{true}}(\bold{x})\right) - 
\hat{f}_{w_{\P_2}}(\bold{x}) \cdot 
\hat{f}_{\mathcal{P}_{2}}(\bold{x})\right|\\
& \quad + \left|\hat{f}_{w_{\P_2}}(\bold{x}) \cdot \hat{f}_{\mathcal{P}_{2}}(\bold{x}) - w_{\mathcal{P}_{2}}(\bold{x}) \cdot \hat{f}_{\mathcal{P}_{2}}(\bold{x})\right|\\
  & \quad +
  \left|
  w_{\P_2}(\bold{x}) \cdot \hat{f}_{\mathcal{P}_{2}}(\bold{x}) - w_{\P_2}(\bold{x}) \cdot f(\bold{x})
  \right|\\
& \leq c_{21} \cdot \left(\max\left\{2a, \|f\|_{C^q(\R^d)}\right\}\right)^{4(q+1)} \cdot \frac{1}{M^{2p}}.
\end{align*}
In case that $\bold{x}$ is contained in \eqref{noset1}, the 
approximation error of $\hat{f}_{\mathcal{P}_{2}}$ is not of size 
$1/M^{2p}$. But the value of $\hat{f}_{\mathrm{check}, 
\mathcal{P}_{2}}(\bold{x})$ is $1$, 
such that  $\hat{f}_{\mathcal{P}_{2}, \mathrm{true}}$ is zero. 
Furthermore, we have 
\begin{align*}
\left|\hat{f}_{w_{\P_2}}(\bold{x})\right| &\leq 2.
\end{align*}
 Thus, $\hat{f}_{w_{\P_2}}(\bold{x})$ and $\hat{f}_{\mathcal{P}_{2}, 
\mathrm{true}}(\bold{x})$ 
 are contained in the interval, where \eqref{appfmult} holds. 
 Together with 
\begin{align*}
0 \leq w_{\P_2}(\bold{x}) \leq \frac{1}{2 \cdot M^{2p}}
\end{align*}
and the triangle inequality it follows 
\begin{align*}
&\left|\hat{f}_{\mathrm{mult}}\left(\hat{f}_{w_{\P_2}}(\bold{x}), 
\hat{f}_{\mathcal{P}_{2}, \mathrm{true}}(\bold{x})\right) - 
w_{\P_2}(\bold{x}) \cdot f(\bold{x})\right|\\
& \leq 
\left|\hat{f}_{\mathrm{mult}}\left(\hat{f}_{w_{\P_2}}(\bold{x}), 
\hat{f}_{\mathcal{P}_{2}, \mathrm{true}}(\bold{x})\right) - 
\hat{f}_{w_{\P_2}}(\bold{x}) \hat{f}_{\mathcal{P}_{2}, 
\mathrm{true}}(\bold{x})\right|\\
& \quad + \left|\hat{f}_{w_{\P_2}}(\bold{x}) \cdot 
\hat{f}_{\mathcal{P}_{2}, \mathrm{true}}(\bold{x}) - 
w_{\P_2}(\bold{x}) \cdot  \hat{f}_{\mathcal{P}_{2}, 
\mathrm{true}}(\bold{x})\right|\\
& \quad + \left|w_{\P_2}(\bold{x}) \cdot  \hat{f}_{\mathcal{P}_{2}, 
\mathrm{true}}(\bold{x}) - w_{\P_2}(\bold{x}) \cdot 
f(\bold{x})\right|\\ 
& \leq c_{22} \cdot \left(\max\{\|f\|_{\infty} , 1\}\right)^2 \cdot \frac{1}{M^{2p}}.
\end{align*}
In case that $\bold{x}$ is in \eqref{noset2}, it is not in \eqref{noset1} 
and the network $\hat{f}_{\mathcal{P}_{2}}$ 
approximates $f(\bold{x})$ with an error as in \eqref{200}. 
Furthermore, $\hat{f}_{w_{\P_2}}(\bold{x}) \in [-2,2]$
approximates $w_{\P_2}(\bold{x})$ with an error as in \eqref{appwv}.
The value of $\hat{f}_{\mathrm{check}, \mathcal{P}_{2}}(\bold{x})$ is 
contained in the interval $[0,1]$, such that
\begin{align*}
\left|\hat{f}_{\mathcal{P}_{2},\mathrm{true}}(\bold{x})\right| \leq 
\left|\hat{f}_{\mathcal{P}_{2}}(\bold{x})\right| \leq 2 \cdot 
\max\left\{\|f\|_{\infty}, 1\right\}.
\end{align*}
Hence $\hat{f}_{w_{\P_2}}(\bold{x})$ and $\hat{f}_{\mathcal{P}_{2}, 
\mathrm{true}}(\bold{x})$ are contained in the interval, where 
\eqref{appfmult} holds. 
Together with 
\begin{align*}
w_{\P_2}(\bold{x}) \leq \frac{1}{M^{2p}}
\end{align*}
and the triangle inequality it follows again
\begin{align*}
&\left|\hat{f}_{\mathrm{mult}}\left(\hat{f}_{w_{\P_2}}(\bold{x}), 
\hat{f}_{\mathcal{P}_{2}, \mathrm{true}}(\bold{x})\right) - 
w_{\P_2}(\bold{x}) \cdot f(\bold{x})\right|\\
& \leq 
\left|\hat{f}_{\mathrm{mult}}\left(\hat{f}_{w_{\P_2}}(\bold{x}), 
\hat{f}_{\mathcal{P}_{2}, \mathrm{true}}(\bold{x})\right) - 
\hat{f}_{w_{\P_2}}(\bold{x}) \cdot \hat{f}_{\mathcal{P}_{2}, 
\mathrm{true}}(\bold{x})\right|\\
& \quad + \left| \hat{f}_{w_{\P_2}}(\bold{x}) \cdot 
\hat{f}_{\mathcal{P}_{2}, \mathrm{true}}(\bold{x}) - 
w_{\P_2}(\bold{x}) \cdot  \hat{f}_{\mathcal{P}_{2}, 
\mathrm{true}}(\bold{x})\right|\\
& \quad + \left|w_{\P_2}(\bold{x}) \cdot  \hat{f}_{\mathcal{P}_{2}, 
\mathrm{true}}(\bold{x})\right|\\
& \leq c_{23} \cdot \left(\max\{\|f\|_{\infty} , 1\}\right)^2 \cdot \frac{1}{M^{2p}}.
\end{align*} 
\hfill $\Box$
\\
\noindent
{\bf Proof of Theorem \ref{th2}.}
Let $\mathcal{P}_1$ and $\mathcal{P}_2$ be the partitions defined as in \eqref{partition}. We set
\begin{align*}
\mathcal{P}_{1,1} = \mathcal{P}_1 \ \text{and} \ \mathcal{P}_{2,1} = 
\mathcal{P}_2
\end{align*}
and define for each $v \in \{2, \dots, 2^d\}$ partitions 
$\mathcal{P}_{1,v}$ 
and $\mathcal{P}_{2,v}$, which are modifications of 
$\mathcal{P}_{1,1}$ and $\mathcal{P}_{2,1}$ where at least one of the 
components it shifted by $1/(2M^2)$.
To avoid that the approximation error of the networks increases close to the boundaries of some cube of the partitions, we multiply each value of $\hat{f}_{\mathcal{P}_{2,v}}$ with a weight 
\begin{align}
\label{w_v}
w_v(\bold{x}) = \prod_{j=1}^d \left(1- 2 \cdot M^2 \cdot 
\left|(C_{\mathcal{P}_{2,v}}(\bold{x}))_{\mathrm{left}}^{(j)} + 
\frac{1}{2 \cdot M^2} - x^{(j)}\right|\right)_+.
\end{align}
It is easy to see that $w_v(\bold{x})$ is a linear tensorproduct B-spline
which takes its maximum value at the center of $C_{\P_{2,v}}(\bold{x})$, which
is nonzero in the inner part of $C_{\P_{2,v}}(\bold{x})$ and which
vanishes
outside of $C_{\P_{2,v}}(\bold{x})$. Consequently
we have $w_1(\bold{x})+ \dots + w_{2^d}(\bold{x}) = 1$
for $\bx \in \R^d$.
Let $\hat{f}_{1}, \dots, \hat{f}_{2^d}$ be the networks of 
Lemma~\ref{le6}
corresponding to the partitions
$\mathcal{P}_{1,v}$ and $\mathcal{P}_{2,v}$ $(v \in \{1, \dots, 2^d\})$, respectively. 
Each $\P_{1,v}$ and $\P_{2,v}$ form a partition of
$\Rd$ and the error bounds of
  Lemma~\ref{le6} hold for each network $\hat{f}_{v}$ on $\M$.
We set
\begin{align*}
\hat{f}_{\mathrm{net}}(\bold{x}) = \sum_{v=1}^{2^d} 
\hat{f}_{v}(\bold{x}).
\end{align*}
Using Lemma~\ref{le6} it is easy to see that this network is contained 
in the network class $\mathcal{F}(L, r)$ with
\begin{align*}
L=5+\lceil \log_4(M^{2p})\rceil \cdot \left(\lceil \log_2(\max\{q,d\} +1)\rceil+1\right) 
\end{align*}
and 
\begin{align*}
  r=2^d \cdot 64 \cdot \binom{d+q}{d} \cdot d^2 \cdot (q+1) \cdot
  \lceil c_{13} \cdot M^{d^*} \rceil.
\end{align*}
Since
\begin{align*}
f(\bold{x}) = \sum_{v=1}^{2^d} w_v(\bold{x}) \cdot f(\bold{x})
\end{align*}
it follows directly by Lemma~\ref{le6}
\begin{align*}
  \left|\hat{f}_{\mathrm{net}}(\bold{x}) - f(\bold{x})\right| 
  &= \left|\sum_{v=1}^{2^d}
  \hat{f}_v (\bx) - 
\sum_{v=1}^{2^d} w_v(\bold{x}) \cdot f(\bold{x})\right|\\
& \leq \sum_{v=1}^{2^d} 
\left|
\hat{f}_v (\bx) -  
w_v(\bold{x}) \cdot f(\bold{x})\right|\\
& \leq c_{24} \cdot \left(\max\left\{3a, \|f\|_{C^q(\R^d)}\right\}\right)^{4(q+1)} \cdot \frac{1}{M^{2p}}.
\end{align*}
\hfill $\Box$
\end{document}